\documentclass[11pt,fleqn]{amsart}

\usepackage{amsmath,amssymb,amsfonts}
\usepackage[round]{natbib}
\usepackage{dsfont}
\usepackage{enumerate}
\usepackage{mathrsfs}
\usepackage{graphicx}
\usepackage{abbreviationsGeneral-ext}
\usepackage{nicefrac}
\usepackage[ansinew]{inputenc}
\usepackage{rotating} % f"ur sidewaystable
\usepackage{array}
\usepackage{dcolumn}

\newcommand{\td}{\tilde{d}}
\newcommand{\tg}{\tilde{g}}

\def\X{\mathbb{X}}
\def\ds{\displaystyle}
\def\NM{N\!M}

\begin{document}

\title{On the efficiency of Gini's mean difference}
\author{Carina Gerstenberger}
\author{Daniel Vogel}
\date{\today}
\address{
Fakult\"at f\"ur Mathematik, Ruhr-Universit\"at Bo\-chum, 
44780 Bochum, Germany}
\email{carina.gerstenberger@rub.de}
\address{
Fakult\"at Statistik, Technische Universit\"at Dortmund, 
44221 Dortmund, Germany}
\email{daniel.vogel@tu-dortmund.de}

\begin{abstract}
\small
The asymptotic relative efficiency of the mean deviation with respect to the standard deviation is 88\% at the normal distribution. 
In his seminal 1960 paper \emph{A survey of sampling from contaminated distributions}, J.~W.~Tukey points out that, if the normal distribution is contaminated by a small $\epsilon$-fraction of a normal distribution with three times the standard deviation, the mean deviation is more efficient than the standard deviation---already for $\epsilon < 1\%$.
This came as a surprise to most statisticians at the time, and the publication is today considered as one of the main pioneering works in the development of robust statistics. 
In the present article, we examine the efficiency of the mean deviation and Gini's mean difference (the mean of all pairwise distances).
% at various distributions. 
The latter is known to have an asymptotic relative efficiency of 98\% at the normal distribution. 
Our findings support the viewpoint that Gini's mean difference combines the advantages of the mean deviation and the standard deviation.  
We also answer the question, what percentage of contamination in Tukey's 1:3 normal mixture model renders Gini's mean difference more efficient than the standard deviation.

\end{abstract}


\maketitle

%\tableofcontents

%  *=======================================================================================*
%  [\       \       \       \       \       \       \       \       \       \       \      ]
%  [ \       \       \       \       \       \       \       \       \       \       \     ]
%  [  \       \       \       \       \       \       \       \       \       \       \    ]
%  [   \       \       \       \       \       \       \       \       \       \       \   ]
%  [    \       \       \       \       \       \       \       \       \       \       \  ]
%  [     \       \       \       \       \       \       \       \       \       \       \ ]
%  [      \       \       \       \       \       \       \       \       \       \       \]
%  [       \       \       \       \       \       \       \       \       \       \       ]
%  *=======================================================================================*
%--------------------------------------------------------------------------------------------------------------------------------------------------------------------------

\section{Introduction}

Let $X$ be a random variable with distribution $F$, and define $F^\star_{a,b}$ as the distribution of $a X + b$.
We call any function $s$ that assigns a non-negative number to any univariate distribution $F$ (potentially restricted to a subset of distributions, e.g.\ with finite second moments) a \emph{measure of variability}, (or a \emph{measure of dispersion} or simply a \emph{scale measure}) if it satisfies
\[
	s(F^\star_{a,b}) = |a|\,s(F)   	
	\qquad \mbox{ for all }  a, b \in \R. 
\]  
In this article, we compare three very common descriptive measures of variability: 
\begin{enumerate}[(i)]
\item the standard deviation $\sigma(F) = \{  E(X-EX)^2 \}^{1/2}$, 
\item the mean absolute deviation (or mean deviation for short) $d(F) = E|X-m(F)|$, where $m(F)$ denotes the median of $F$, and
\item Gini's mean difference $g = E|X-Y|$.
\end{enumerate}
Here, $X$ and $Y$ are independent and identically distributed random variables with distribution function $F$.  Recall that the variance can also be written as $\sigma^2(F) = E(X-Y)/2$. 
%We identify the distribution with its cdf and do not make any distinction in the notation.
We define the median $m(F)$ as the center point of the set $\{ x \in \R \, |\, F(x-) \le 1/2 \le F(x)\}$, where $F(x-)$ denotes the left-hand side limit.
%We will call the mean absolute deviation simply the mean deviation in the following, the abbreviation MAD commonly refers to the median absolute deviation \citep{Hampel:1974}. 
%Gini's mean difference is commonly attributed to \citet{Gini...}
%
%
%
%
Suppose now we observe data $\X_n = (X_1,\ldots,X_n)$, where the $X_i$, $i  = 1, \ldots, n$, are independent and identically distributed with cdf $F$. Let $\hat{F}_n$ be the corresponding empirical distribution function. The natural estimates for the above scale measures are the functionals applied to $\hat{F}_n$. However, we define the sample versions of the standard deviation and the mean deviation slightly differently. Let
\begin{enumerate}[(i)]
\item 
$\displaystyle \sigma_n = \sigma_n(\X_n) = \Big\{ \frac{1}{n-1} \sum_{i=1}^n \left( X_i - \bar{X}_n \right)^2 \Big\}^{1/2}  
																					= \Big\{ \frac{1}{n(n-1)} \sum_{1\le i < j \le n} \left(X_i-X_j\right)^2 \Big\}^{1/2} $ \\
denote the sample standard deviation, 
\item 
$\displaystyle d_n = d_n(\X_n) = \frac{1}{n-1} \sum_{i=1}^n |X_i - m(\hat{F}_n)|$  \ the sample mean deviation and 
\item 
$\displaystyle g_n  = g_n(\X_n) = \frac{2}{n(n-1)} \sum_{1 \le i < j \le n} |X_i - X_j|$ \ the sample mean difference. 
\end{enumerate}

While it is common practice to use $1/(n-1)$ instead of $1/n$ in the definition of the sample variance, due to the thus obtained unbiasedness, it is not so clear what finite-sample version of the mean deviation to use. Unfortunately, the factor $1/(n-1)$ does not yield unbiasedness for any distribution, as it is the case for the variance, but it leads to a significantly smaller bias in all our finite-sample simulations, see Section \ref{sec:finite-sample}. 

Furthermore, there is the question of the location estimator, which applies, in principle, to the mean deviation as well as to the standard deviation, and also to their population versions. While it is again established to use the mean along with the standard deviation, the picture is less clear for the mean deviation. We propose to use the median, mainly due to conceptual reasons: the median minimizes the mean deviation as the mean minimizes the standard deviation. This also suggests to apply the simple $1/(n-1)$ bias correction in both cases. However, our main results concern asymptotic efficiencies at symmetric distributions, for which the choice of the location measure as well as $n$ vs.\ $n-1$ question is irrelevant. 
%We also found little differences in the finite-sample simulations. 
%
%
%

If $E X^2 < \infty$, Gini's mean difference and the mean deviation are asymptotically normal. For the asymptotic normality of $\sigma_n$, fourth moments are required. Strong consistency and asymptotic normality of $g_n$ and $\sigma_n^2$ follow from general $U$-statistic theory \citep{hoeffding:1948}, and thus for $\sigma_n$ by a subsequent application of the continuous mapping theorem and the delta method, respectively.
Letting
\[
	d_n(\X_n,t) = \frac{1}{n-1}\sum_{i=1}^n |X_i - t|,
\]
the asymptotic normality of $d_n(\X,t)$ for any fixed location $t$ holds also under the existence of second moments and is a simple corollary of the central limit theorem. The asymptotic normality of $d_n(\X_n,t_n)$, where $t_n$ is a location estimator is not equally straightforward \citep[cf.\ e.g.][Theorem 5 and the examples below]{Bickel1975}. A set of sufficient conditions is that $\sqrt{n}(t_n-t)$ is asymptotically normal and $F$ is symmetric around $t$.

%% Wozu das alles %%%%
%This article is a contribution to the area of robust statistics, but it is about non-robust estimators. 
The standard deviation is, with good cause, the by far most popular measure of variability. One main reason for considering alternatives is its lack of robustness, i.e.\ its susceptibility to outliers and its low efficiency at heavy-tailed distributions. 
The two alternatives considered here are --- in the modern understanding of the term --- not robust, but they are more robust than the standard deviation. The extreme non-robustness of 
%the sample variance and 
the standard deviation, which also emerges when comparing it with the mean deviation, played a vital role in recognizing the need for robustness and thus helped to spark the development of robust statistics, cf.\ e.g.\ \citet{Tukey1960}. 
The purpose of this article is to introduce Gini's mean difference into the old debate of mean deviation vs.\ standard deviation \citep[e.g.][]{Gorard2005} --- not as a compromise, but as a consensus. We will argue that Gini's mean difference combines the advantages of the standard deviation and the mean deviation. 

%%%% Was passiert genau und was kommt raus
When proposing robust alternatives to any normality-based standard estimator, 
%(such as the standard deviation here), 
the gain in robustness is usually paid by a loss in efficiency at the normal model. The two aspects, robustness and efficiency, have to be analyzed and be put into relation with each other. 
The theoretical robustness properties of the three estimators are quickly summarized: they all have an asymptotic breakdown point of zero and an unbounded influence function. There are some slight advantages for the mean deviation and Gini's mean difference: their influnce functions increase linearly (as compared to the quadratic increase for the standard deviation), and they require only second moments to be asymptotically normal (as compared to the 4th moments for the standard deviation). The influence functions of all three estimators at the standard normal distribution are plotted in Figure \ref{fig:if}.

We are thus left to study their efficiencies. This is the main concern in this paper. We compute and compare the asymptotic variances of the estimates at several distributions. We restrict our attention to symmetric distributions, since we are interested primarily in the effect of the tails of the distribution, which arguably have the most decisive influence on the behavior of the estimators. 
We consider in particular the $t_\nu$ distribution and the normal mixture distribution, which are both popular outlier models in robust statistics. To summarize our findings, in all relevant situations where Gini's mean difference does not rank first among the three estimators in terms of efficiency, it does rank second with very little difference to the respective winner. A more detailed discussion is deferred to Section \ref{sec:conclusion}.

The rest of the paper is organized as follows: In Section \ref{sec:asymptotic}, asymptotic efficiencies of the scale estimators are compared. We study in particular their asymptotic variances at the normal mixture model. In Section \ref{sec:if}, the influence functions are computed.
We complement our findings with finite-sample simulations in Section \ref{sec:finite-sample}. Section \ref{sec:conclusion} contains a summary. 
Remarks on the computation of the asymptotic variances are given in the Appendix.

We close this section by introducing some further terms and notation. 
Letting $s_n$ be any of the estimators above and $s$ the corresponding population value, we define the asymptotic variance $ASV(s_n) = ASV(s_n;F)$ of  $s_n$ at the distribution $F$ as the variance of the limiting normal distribution of $\sqrt{n}(s_n- s)$, when $s_n$ is evaluated at an independent sample $X_1,\ldots,X_n$ drawn from $F$. We note that, in general, convergence in distribution does not imply convergence of the second moments without further assumptions (uniform integrability), but it is usually the case in situations encountered in statistical applications, specifically it is true for the estimators considered here, and we may write 
\[
		ASV(s_n) = \lim_{n\to\infty} n\, \var(s_n).
\]
We are going to compute asymptotic relative efficiencies of $g_n$ and $d_n$ with respect to $\sigma_n$. Generally, for two estimators $a_n$ and $b_n$ with $a_n \cip \mu \quad \mbox{and} \quad b_n \cip \mu$ for some $\mu \in \R$, the asymptotic relative efficiency of $a_n$ with respect to $b_n$ at distribution $F$ is defined as 
\[
	ARE(a_n,b_n; F) = ASV(b_n;F)/ASV(a_n;F).
\]
In order to make the scale estimators comparable efficiency-wise, we introduce the $\sigma$-standardized versions of the estimators, 
\[
	\tg_n = \frac{\sigma}{g} g_n, 
	\qquad \td_n = \frac{\sigma}{d} d_n.
%	\qquad \mbox{and} \quad
%	\qquad \tr_n = \frac{\sigma}{r} r_n
\]
These estimators are of no practical use, since they require the knowledge of the parameter $\sigma$, which they aim to estimate, but since they  estimate the same quantity as $\sigma_n$, we may compare their asymptotic variances. We then define the asymptotic relative efficiency of $s_n$ (where $s_n$ may be any scale estimator) with respect to the standard deviation at the population distribution $F$ as
\be \label{eq:are}
	ARE(s_n,\sigma_n;F) \ = \ ARE(\tilde{s}_n,\sigma_n;F) \ = \ \frac{ASV(\sigma_n;F)}{ASV(s_n;F)} \frac{s^2(F)}{\sigma^2(F)}. 
\ee
Also, if we compare the efficiencies of two scale estimators $s_n^{(1)}$ and $s_n^{(2)}$, the comparison shall refer to their $\sigma$-standardized versions.	

%
%
%  *=======================================================================================*
%  [\       \       \       \       \       \       \       \       \       \       \      ]
%  [ \       \       \       \       \       \       \       \       \       \       \     ]
%  [  \       \       \       \       \       \       \       \       \       \       \    ]
%  [   \       \       \       \       \       \       \       \       \       \       \   ]
%  [    \       \       \       \       \       \       \       \       \       \       \  ]
%  [     \       \       \       \       \       \       \       \       \       \       \ ]
%  [      \       \       \       \       \       \       \       \       \       \       \]
%  [       \       \       \       \       \       \       \       \       \       \       ]
%  *=======================================================================================*

\section{Asymptotic efficiencies}
\label{sec:asymptotic}

We gather the general expressions for the population values and asymptotic variances of the three scale measures (Section \ref{subsec:gen}) and then evaluate them at several symmetric distributions (Section \ref{subsec:spec}). 
We study the two-parameter family of the normal mixture model in some detail in Section \ref{subsec:nm}.

%  *=======================================================================================*
%  [\       \       \       \       \       \       \       \       \       \       \      ]
%  [ \       \       \       \       \       \       \       \       \       \       \     ]
%  [  \       \       \       \       \       \       \       \       \       \       \    ]
%  *=======================================================================================*

\subsection{General expressions}
\label{subsec:gen}

The exact finite-sample variance of the empirical variance $\sigma^2_n$ is 
\[
	\var(\sigma_n^2) \ = \ \frac{1}{n}\left\{ \mu_4 - 4\mu_3 \mu_1 + 3\mu_2^2 - 2\sigma^2 \frac{2n-3}{n-1} \right\},
\]
where $\mu_k = EX^k$, $k \in \N$, is the $k$th non-central moment of $X$, in particular $\sigma^2 = \sigma^2(F) = \mu_2 - \mu_1^2$. Thus
$ASV(\sigma^2_n)	= \mu_4  + 3\mu_2^2 - 4\left\{ \mu_3 \mu_1 + \sigma^4\right\}$, and hence we have by the delta method
\be \label{eq:ASV.sigma_n}
	ASV(\sigma_n) = \frac{\mu_4   - 4\mu_3 \mu_1 + 3\mu_2^2}{4\sigma^2} - \sigma^2.
\ee
If the distribution $F$ is symmetric around $E(X)=\mu_1$ and has a Lebesgue density $f$, the mean deviation $d = d(F)$ can be written as
\be \label{eq:d.sym}
	d \ = \int_{-\infty}^{\infty} |x-\mu_1| f(x)\,  dx
	 \ = \ 2 \int_{\mu_1}^{\infty} (x-\mu_1) f(x)\,  dx
\ee
The asymptotic variances of $d_n$ and its $\sigma$-standardized version $\td_n$ are 
$ASV(d_n) = \sigma^2 - d^2$ and $ASV(\td_n) = \sigma^2 \left\{ \sigma^2/d^2 - 1 \right\}$, respectively. 
For any $F$ possessing a Lebesgue density $f$, Gini's mean difference $g = g(F)$ is
\be \label{eq:g}
	  g \ =  \ \int_{-\infty}^{\infty} \int_{-\infty}^{\infty} |x - y|\, f(x)\, f(y)\, dy\, dx 
	    \ =  \ 2 \int_{-\infty}^{\infty} \int_{x}^{\infty} (y - x)\, f(x)\, f(y)\,  dy\, dx, 
\ee
which can be further reduced to  
\be \label{eq:g.sym}
	g \ = \ 4 \int_{-\infty}^{\infty} \int_{x}^{\infty} y \, f(y) \, dy\, f(x)\,  dx
	 \ = \ 8 \int_{0}^{\infty} \int_{x}^{\infty} y \, f(y) \, dy\, f(x)\,  dx
\ee
if $F$ is symmetric around 0.
\citet{lomnicki:1952} gives the variance of the sample mean difference $g_n$ as
\be \label{eq:var.g_n}
	\var(g_n) = \frac{1}{ n(n-1) }
		\left\{ 
			4(n-1) \sigma^2 + 16 (n-2) J - 2 (2n-3) g^2 
		\right\},
\ee
where
\be \label{eq:J}
	J = \int_{x = -\infty}^{\infty} \int_{y = -\infty}^{x} \int_{z = x}^{\infty} 
			(x-y)(z-x) f(z) f(y) f(x)\, dz\, dy\, dx. 
\ee
Thus, the asymptotic variances of $g_n$ and its $\sigma$-standardized version $\tg_n$ are 
$	ASV(g_n) = 4 \{ \sigma^2 + 4 J - g^2  \}$ and 	$ASV(\tg_n) = 4 \sigma^2\{ (\sigma^2 + 4J)/g^2 - 1 \}$, 
respectively.

%
%
%
%
%
%
%
%
%
%
%
%
%  *=======================================================================================*
%  [\       \       \       \       \       \       \       \       \       \       \      ]
%  [ \       \       \       \       \       \       \       \       \       \       \     ]
%  [  \       \       \       \       \       \       \       \       \       \       \    ]
%  *=======================================================================================*
%
\subsection{Specific distributions}
\label{subsec:spec}

Table~\ref{tab:specific.1} lists the densities and first four moments of the following distribution families: normal, Laplace, uniform, $t_\nu$ and normal mixture. The resulting expressions for $\sigma(F)$, $d(F)$ and the asymptotic variances of their sample versions are given in Table~\ref{tab:specific.2}, and for Gini's mean difference, including the integral $J$, in Table~\ref{tab:specific.3}. 
While the contents of Table \ref{tab:specific.2} are straightforward and stated here without proof, the results for Gini's mean difference require the evaluation of the integrals (\ref{eq:g.sym}) and (\ref{eq:J}), which is non-trivial for many distributions. Details for the $t_\nu$ and the normal mixture distribution are given in the Appendix. The expressions for the normal case are due to \citet{nair:1936}. 
For convenience, resulting numerical values of the three scale measures and their asymptotic variances are listed in Table~\ref{tab:numeric.1}. 
Table~\ref{tab:numeric.2} contains the asymptotic relative efficiencies, cf.~(\ref{eq:are}). In particular, we have at the normal model
\[
	ARE(g_n, \sigma_n) = \left\{ \frac{2}{3} \pi + 4 (\sqrt{3}-2) \right\}^{-1} = 0.9779, 
	\qquad 
	ARE(d_n, \sigma_n) = \frac{1}{\pi - 2} = 0.876,
\]
and at the Laplace (or double exponential) distribution
\[
	ARE(g_n, \sigma_n) = 135/112 = 1.2054,  	
	\qquad 
	ARE(d_n, \sigma_n) = 5/4.
\]
Thus, in both situations, Gini's mean difference has an efficiency of more than 96\% with respect to the respective maximum likelihood estimator.
Furthermore, we observe that Gini's mean difference $g_n$ is asymptotically more efficient than the standard deviation $\sigma_n$ at $t_\nu$ distribution for $\nu \le 40$. The mean deviation $d_n$ is asymptotically more efficient than $\sigma_n$ for $\nu \le 15$ and more efficient than $g_n$ for $\nu \le 8$. Thus in the range $9 \le \nu \le 40$, Gini's mean difference is the most efficient of the three. 

One can view the uniform distribution as a limiting case of very light tails. While our focus is on heavy-tailed scenarios, we include the uniform distribution in our study as a simple approach to compare the estimators under light tails. We find a similar picture as under normality: Gini's mean difference and the standard deviation perform equally well, while the mean deviation has a substantially lower efficiency. 
However, it must be noted that the uniform distribution itself is rarely encountered in practice. The limited range is a very strong information, which allows a super-efficient inference.% for the parameters of this family. 

We also include the interquartile range (the distance between the upper and the lower quartile) in the efficiency comparison of Table~\ref{tab:numeric.2} --- without examining this estimator in detail. The purpose is to give a rough impression of how the numbers given compare to another well-known scale measure.  This comparison, though, must be into perspective with two aspects: Firstly, the interquartile range is primarily used as a descriptive statistic for data sets rather than an estimator for a true population value. Secondly, the interquartile range is much more robust, it has a bounded influence function and a breakdown point of about 0.25. Also, there are other highly robust scale measures which are more efficient than the interquartile range, for instance the median absolute deviation \citep[MAD,][]{Hampel1974} or the $Q_n$ by \citet{RousseeuwCroux1993}. We do not attempt to give a complete review, which is clearly beyond the scope of this paper.

Finally, we take a closer look at the normal mixture distribution and explain our choices for $\lambda$ and $\epsilon$ in Table~\ref{tab:numeric.2}.
\begin{table}
\caption{
		Densities and non-central moments of several parametric families. The scaling factor for the $t_\nu$ distribution is
		$c_\nu = \Gamma (\frac{\nu+1}{2})/(\sqrt{\nu \pi}\,\Gamma( \frac{\nu}{2}))$.
}
\newcolumntype{C}[1]{>{\centering}m{#1}}
%----------------
\newdimen\mylength
\mylength=0.33\textwidth
%---------------
% wird unten 
\renewcommand{\arraystretch}{2.5}
\begin{center}
\begin{tabular}{C{0.13\textwidth}|C{0.25\textwidth}|C{0.18\textwidth}|c}
%\begin{tabular}{c|c|c|c}
distribution & density $f(x)$ & parameters & moments \\
\hline
% Normalverteilung
		normal & 
		$\frac{1}{\sqrt{2\pi \sigma^2} } \exp\left\{ - \frac{(x-\mu)^2}{2 \sigma^2} \right\}$ &
		$\mu \in \R, \sigma^2 > 0$  & 
		\begin{minipage}[c]{\mylength}
			\renewcommand{\baselinestretch}{1.2} \normalsize
			\centering
			\smallskip
			$\mu_1 = \mu$, \ $\mu_2 = \sigma^2 + \mu^2$, 
			
		%	\smallskip
			$\mu_3 = \mu^3 + 3 \mu \sigma^2$, 
			
			$\mu_4 = \mu^4 + 6 \mu^2 \sigma^2 + 3\sigma^4$ 
		\end{minipage} \\[3.0ex]
\hline		
% Laplace-Verteilung
		Laplace & 	
		$\frac{1}{2\alpha}\exp\left\{\frac{-|x-\mu|}{\alpha}\right\}$ &
		$\mu \in \R, \alpha >0$  & 
		\begin{minipage}{\mylength}
			\smallskip
			\renewcommand{\baselinestretch}{1.2} \normalsize
			\centering
			$\mu_1 = \mu$, \ $\mu_2 = \mu^2+2\alpha^2$, 
		
			$\mu_3 = \mu^3+6\alpha^2\mu$, 
		
			$\mu_4 = \mu^4+12\alpha^2\mu^2+24\alpha^4$ 
		\end{minipage} \\[3.0ex]
\hline		
% Uniform-Verteilung
		uniform & 
		$\ds \frac{1}{b-a} \Varind{[a,b]}(x)$ &
		{\small $-\infty < a < b < \infty$} & 
		\begin{minipage}{\mylength}
			\smallskip
			\renewcommand{\baselinestretch}{1.3} \normalsize
			\centering
			$\mu_1 = \frac{1}{2}(a+b)$,  
			
			$\mu_2 = \frac{1}{3}\left\{(a+b)^2-ab\right\}$, 
			
			$\mu_3 = \frac{1}{4}(a+b)(a^2+b^2)$, 
			
			$\mu_4 = \frac{1}{5}\left\{(a+b)(a^3+ab^2)+b^4\right\}$  	
		\end{minipage}\\[5.5ex]	
\hline
% t-Verteilung
		$t_\nu$ & 	
		$\ds c_\nu \left( 1 + \frac{x^2}{\nu} \right)^{-\frac{\nu+1}{2}}$  &
		$\nu \in \N$  & 
		\begin{minipage}{\mylength}
			\smallskip
			\renewcommand{\baselinestretch}{1.2} \normalsize
			\centering
			$\mu_1 = \mu_3 = 0$, 
			
			$\mu_2 = \nu/(\nu-2)$, 
			
			$\ds \mu_4 = 3 \nu^2/\{(\nu-2)(\nu-4)\} $ %\frac{3 \nu^2}{(\nu-2)(\nu-4)}$ 
		\end{minipage} \\[3.0ex]		
\hline	 
% Normal-Mixture-Verteilung
		normal mixture & 
		$\epsilon \frac{1}{\sqrt{2\pi}\lambda}\exp{\{-\frac{x^2}{2\lambda^2}\}}+\left(1-\epsilon\right)\frac{1}{\sqrt{2\pi}}\exp{\{-\frac{x^2}{2}\}}$ &
		$0 \le \epsilon \le 1$, $\lambda \ge 1$ & 
		\begin{minipage}{\mylength}
			\smallskip
			\renewcommand{\baselinestretch}{1.2} \normalsize
			\centering
			$\mu_1 = \mu_3 = 0$, 
			
			$\mu_2 = \epsilon\lambda^2+\left(1-\epsilon\right)$, 
			
			$\mu_4 = 3\epsilon\lambda^4+3\left(1-\epsilon\right)$ 
		\end{minipage} \\[2.0ex]		
%\hline	 
\end{tabular}
\end{center}
\label{tab:specific.1}
\end{table}
%
%
%

%\vfill
\begin{table}
%\medskip
\bigskip
\caption{Specific values of $\sigma$, $d$ and the respective asymptotic variances for the distribution families given in Table \ref{tab:specific.1}. $c_\nu = \Gamma (\frac{\nu+1}{2})/(\sqrt{\nu \pi}\,\Gamma( \frac{\nu}{2}))$}
\newcolumntype{C}[1]{>{\centering}m{#1}}
\renewcommand{\arraystretch}{2.5}
%\begin{tabular}{C{0.13\textwidth}|C{0.2\textwidth}|C{0.25\textwidth}|C{0.18\textwidth}|c}
\begin{center}
\begin{tabular}{C{0.12\textwidth}|c|c|c|c}
	distribution  & $\sigma(F)$ & $ASV(\sigma_n$) & $d(F)$ &  $ASV(d_n)$ \\
\hline
% Normalverteilung
		normal & 
		$\sigma$ & 
		$\ds \frac{\sigma^2}{2}$ & 
		$\ds\frac{2\sigma}{\sqrt{2\pi}}$  &  
		$\ds \sigma^2 \left\{ 1 - \frac{2}{\pi} \right\}$ \\%[2.0ex]
\hline		
% Laplace-Verteilung
		Laplace  & 
		$\ds\sqrt{2}\alpha$ & 
		$\ds \frac{5}{2}\alpha^2$ & 
		$\alpha$ 		& 
		$\alpha^2$ \\%[2.0ex]
\hline		
% Uniform-Verteilung
		uniform & 
		$\ds \frac{b-a}{2\sqrt{3}}$ & 
		$\frac{1}{60}(b-a)^2$ & 
		$\ds \frac{b-a}{4}$ & 
		$\frac{1}{48}(b-a)^2$ \\%[2.0ex]		
\hline
% t-Verteilung
		$t_\nu$ & 
		$\ds \sqrt{\frac{\nu}{\nu-2}}$ & 
		$\ds \frac{ \nu (\nu-1)}{ 2 (\nu-2)(\nu-4)}$ & 
		$\ds \frac{2\nu c_\nu}{\nu - 1}$ & 
		$\ds \frac{\nu}{\nu - 2} - \left\{ \frac{2\nu c_\nu}{\nu - 1} \right\}^2$ \\%[2.0ex]		
\hline	 
% Normal-Mixture-Verteilung
		normal mixture & 
		$\sqrt{\epsilon\lambda^2+\left(1-\epsilon\right)}$ & 
    $\frac{3\left(\epsilon\lambda^4+1-\epsilon\right)-\left(\epsilon\lambda^2+1-\epsilon\right)^2}{4\left(\epsilon\lambda^2+1-\epsilon\right)}$ & 
		$\sqrt{\frac{2}{\pi}} \{ \epsilon\lambda + \left(1-\epsilon\right)\}$ & 
			\begin{minipage}{0.20\textwidth}
			\smallskip
			\renewcommand{\baselinestretch}{1.2} \normalsize
			\centering
			$\epsilon\lambda^2 + 1 - \epsilon$
			
			$-\frac{2}{\pi}\{ \epsilon\lambda+\left(1-\epsilon\right)\}^2$ 
		\end{minipage}\\%[2.0ex]		
%\hline	 
\end{tabular}
\end{center}
\label{tab:specific.2}
\end{table}
\begin{table}
\caption{Population values, cf.~(\ref{eq:g}), expressions for $J$, cf.~(\ref{eq:J}), and resulting asymptotic variances for Gini's mean difference at the parametric families of Table~\ref{tab:specific.1}. 
Abbreviations used:		$\zeta(\lambda) = \sqrt{2+\lambda^2}$,
  $K_\nu \ = \int_{-\infty}^{\infty} x^2 f_\nu(x) F_\nu^2(x) \, dx$ with $f_\nu$, $F_\nu$ being density and cdf of the $t_\nu$ distribution.
  $B(\cdot,\cdot)$ denotes the beta function.}
\newcolumntype{C}[1]{>{\centering}m{#1}}
\renewcommand{\arraystretch}{2.5}
%\begin{tabular}{C{0.13\textwidth}|C{0.28\textwidth}|C{0.38\textwidth}|c}
\begin{center}
\begin{tabular}{C{0.12\textwidth}|c|c|c}
%\hline
distribution & $g(F)$ &   $J$   & $ASV(g_n)$ \\
\hline
%------------------------------------------------------------------------------------------------------------------------------------
% Normal Verteilung
	normal
	&  $\ds \frac{2 \sigma}{\sqrt{\pi}}$ 
	&  $\ds \left(\frac{\sqrt{3}}{2\pi} - \frac{1}{6}\right)\sigma^2$
	&  {\small $ \big\{  \frac{4}{3} + \frac{8}{\pi}(\sqrt{3}-2) \big\}\sigma^2$} \\[2.0ex]
\hline
%------------------------------------------------------------------------------------------------------------------------------------
% Laplace Verteilung
	Laplace 
	& $\ds \frac{3}{2}\alpha$ 
	& $\ds \frac{5}{24}\alpha^2$ 
	& $\ds \frac{7}{3}\alpha^2$ \\[2.0ex]
\hline
%------------------------------------------------------------------------------------------------------------------------------------
% Uniform Verteilung
	uniform 
	& $\ds \frac{1}{3}(b-a)$ 
	& $\ds \frac{1}{120} (b-a)^2$ 
	&  $\ds \frac{1}{45} (b-a)^2$ \\[2.0ex]
\hline
%------------------------------------------------------------------------------------------------------------------------------------
% t Verteilung
	$t_\nu$ 
%	& {\small $\! \frac{4\sqrt{\nu/\pi} }{\nu-1}  
%		\left\{ \! \frac{ \Gamma(\frac{\nu+1}{2}) }{ \Gamma(\frac{\nu}{2})} \!\right\}^{\!2}
%		\frac{ \Gamma\left(\nu-\frac{1}{2}\right) }{ \Gamma\left(\nu \right) } \! $ }
	& $\frac{4\sqrt{\nu}}{\nu-1} 
		\frac{B\big(\frac{\nu}{2}+\frac{1}{2},\, \nu - \frac{1}{2}\big)}{ B\big(\frac{\nu}{2}, \frac{1}{2}\big) B\big(\frac{\nu}{2}, \,\nu\big)   }$
	& $ \frac{2\, \nu}{(\nu-1)^2}
	   \frac{ B\big(\frac{3 \nu}{2}-1,\, \frac{1}{2}\big) }{ B\big(\frac{\nu}{2}, \frac{1}{2}\big)^3}
	   -\frac{\nu}{2(\nu-2)}+K_\nu$
	& {\small $4 \{ \sigma^2 \!+\! 4 J \!-\! g^2 \}$} \\[2.0ex]
\hline
%------------------------------------------------------------------------------------------------------------------------------------
% normal mixture Verteilung
	normal mixture 
	&   \begin{minipage}{0.2\textwidth}
			 \smallskip
			 \renewcommand{\baselinestretch}{1.2} \small
			 \centering
			 			$\ds   \frac{2}{\sqrt{\pi}} \Big\{ \lambda\epsilon^2 + (1-\epsilon)^2 \, + $
			 			
	   				$ \ds \epsilon(1-\epsilon)\sqrt{2\left(1+\lambda^2\right)} \Big\}$
	   	\end{minipage}
	   %
	   %
	   % das ominöse J!
  &  \begin{minipage}{0.38\textwidth}
			 \smallskip
			 \renewcommand{\baselinestretch}{1.2} \small
			 %\centering
       %\footnotesize 
    $\ds (\frac{1}{3}+\frac{\sqrt{3}}{2\pi}) \{ \epsilon^3 \lambda^2 + (1-\epsilon)^3 \}  - \frac{\epsilon \lambda^2 + 1 - \epsilon}{2}$
    
    \smallskip
    $\ds \quad + \, \epsilon^2(1-\epsilon)\bigg[ \frac{\lambda^2}{2} + \frac{1}{4} + \frac{3 \lambda \zeta(\lambda)}{2 \pi} $ \hfill {}

    %$\ds \qquad \quad + \, \frac{1}{\pi \zeta(\lambda) \xi(\lambda)} \Big\{ \frac{(\lambda^2 +3)\lambda^3}{2} + \lambda \Big\} $
    
    \smallskip
    $ \qquad \quad  + \, \frac{\lambda^2}{\pi} \atan(\frac{\lambda}{\zeta(\lambda)}) + \frac{1}{2\pi} \atan( \frac{1}{\lambda \zeta(\lambda)} )\bigg] $ 

    \smallskip 
    $\ds \quad + \, \epsilon(1-\epsilon)^2\bigg[ \frac{\lambda^2}{4} + \frac{1}{2} + \frac{3 \sqrt{1+2\lambda^2}}{2\pi}$
    
    %\smallskip
    %$\ds \qquad \quad + \,\frac{1}{\pi \zeta(1/\lambda) \xi(1/\lambda)} \Big\{\frac{(\lambda^{-2} +3)}{2\lambda} + \lambda\Big\}$
     
     \smallskip
    $ \qquad \ + \, \frac{\lambda^2}{2\pi}  \atan(\frac{\lambda}{\zeta(1/\lambda)}) + \frac{1}{\pi} \atan( \frac{1}{\lambda \zeta(1/\lambda)} )  \bigg]$      	
    \smallskip
    \end{minipage}
		& {\small $4 \{ \sigma^2 \!+\! 4 J \!-\! g^2 \}$ } \\%[12.0ex]
\hline
\end{tabular}
\end{center}
\label{tab:specific.3}
\end{table}

\begin{table}[ht]
\renewcommand{\arraystretch}{1.2}
\centering
\caption{Values and asymptotic variances of the standard deviation $\sigma$, Gini's mean difference $g$ and the mean absolute deviation $d$ at the standard normal distribution $N(0,1)$, the standard Laplace distribution $L(0,1)$, the uniform distribution $U(0,1)$ and several members of the $t_\nu$ family and the normal mixture family $\NM(\lambda,\epsilon)$.}
%\medskip
\begin{tabular}{c| D{.}{.}{6}@{\quad} D{.}{.}{6}@{\quad} D{.}{.}{6}| D{.}{.}{6}@{\quad} D{.}{.}{6}@{\quad} D{.}{.}{6}} 
%-----------------------------------------------------------------------------------
distribution  															& 
\multicolumn{1}{c}{$\sigma$}								& 
\multicolumn{1}{c}{$g$}				 							& 
\multicolumn{1}{c|}{$d$ }										&
\multicolumn{1}{c}{\small $ASV(\sigma_n)$} 	&
\multicolumn{1}{c}{\small $ASV(g_n)$} 			&
\multicolumn{1}{c}{\small $ASV(d_n)$}			\\
\hline
%-----------------------------------------------------------------------------------
	$N(0,1)$              & 1					& 1.128379 		& 0.797884 	          & 0.5				& 0.651006 		& 0.36338 \\
	$L(0,1)$             	& 1.414214	& 1.5				  & 1				           	& 2.5				& 2.333333		& 1       \\
	$U(0,1)$             	& 0.288675	&	0.333333		& 0.25			         	& 0.016667	&	0.022222		& 0.020833  \\
\hline
	$t_5$                 & 1.290994  & 1.383983  	& 0.949017      & 3.333333  & 1.784415  	& 0.766034 \\
	$t_6$                	& 1.224745  & 1.331554	 	& 0.918559 	   	& 1.875		  & 1.453316	 	& 0.656250 \\
	$t_7$                	& 1.183216 	& 1.29694		 	& 0.898313   		& 1.4			 	& 1.268881	 	& 0.593033 \\ 
  $t_{10}$            	& 1.118034 	& 1.239891	 	& 0.864685   		& 0.9375	 	& 1.013824	 	& 0.502319 \\ 
 	$t_{15}$            	& 1.074172 	& 1.199657	 	& 0.840757   		& 0.734266 	& 0.864755	 	& 0.446974 \\
\hline
 	$t_{16}$            	& 1.069045 	& 1.194859	 	& 0.837891 	  	& 0.714286 	& 0.848464	 	& 0.440796 \\
 	$t_{25}$            	& 1.042572 	& 1.169776	 	& 0.822862	   	& 0.621118 	& 0.768006	 	& 0.409855 \\
  $t_{40}$            	& 1.025978 	& 1.153794	 	& 0.813245 	  	& 0.570175 	& 0.720625	 	& 0.391264 \\
 	$t_{41}$            	& 1.025320 	& 1.153156	 	& 0.812861 	  	& 0.568261 	& 0.718794	 	& 0.390540 \\
	$t_{100}$            	& 1.010153 	& 1.138367	 	& 0.803932      & 0.526148 	& 0.677563	 	& 0.374102 \\
\hline
%	$\NM(1.5,0.008)$        & 1.004988  & 1.133336    & 0.801076     & 0.514208  & 0.663118    & 0.368277 \\
%	$\NM(1.5,0.00175)$      & 1.001093  & 1.129464    & 0.798583     & 0.503136  & 0.653665    & 0.364453 \\
%	$\NM(1.5,0.000309)$     & 1.000193  & 1.128571    & 0.798008     & 0.500555  & 0.651476    & 0.363570 \\
	$\NM(3,0.008)$          & 1.031504  & 1.150661    & 0.810651     & 0.890015  & 0.791360     & 0.406845 \\
	$\NM(3,0.00175)$        & 1.006976  & 1.133259    & 0.800677     & 0.589695  & 0.681919    & 0.372916 \\
	$\NM(3,0.000309)$       & 1.001235  & 1.129241    & 0.798378     & 0.516028  & 0.656474    & 0.365065 \\
\end{tabular} 
\label{tab:numeric.1}
\end{table}
\begin{table}[ht]
\caption{Asymptotic relative efficiencies of Gini's mean difference $g_n$, the mean deviation and the interquartile range ($IQR$) with respect to the standard deviation at the standard normal distribution $N(0,1)$, the standard Laplace distribution $L(0,1)$, the uniform distribution $U(0,1)$ and several members of the $t_\nu$ family and the normal mixture family $\NM(\lambda,\epsilon)$.}
\renewcommand{\arraystretch}{1.2}
\centering
\begin{tabular}{c|@{\quad } D{.}{.}{6}@{\qquad } D{.}{.}{6}@{\qquad } D{.}{.}{6}@{\qquad } D{.}{.}{6}}  %-----------------------------------------------------------------------------------
distribution  		&
\multicolumn{1}{c}{\small $ARE(g_n,\sigma_n)$} &
\multicolumn{1}{c}{\small $ARE(d_n,\sigma_n)$} &
\multicolumn{1}{c}{\small $ARE(IQR,\sigma_n)$}	 \\
\hline
%-----------------------------------------------------------------------------------
	$N(0,1)$            & 0.977901		& 0.875969 		& 0.367529 			\\
	$L(0,1)$            & 1.205357	& 1.25				& 0.600566			\\
	$U(0,1)$            & 1					&	0.6					& 0.2						\\
\hline
	$t_5$               & 2.146820  & 2.351417  	& 1.333182  		\\
	$t_6$ 		          & 1.524991  & 1.607143	 	& 0.847705 			\\
	$t_7$				        & 1.325620	& 1.360745	 	& 0.686402 		 	\\
  $t_{10}$          	& 1.137276 	& 1.116343	 	& 0.525866 	 		\\
 	$t_{15}$          	& 1.059075 	& 1.006384	 	& 0.453437 	 		\\
\hline 	
 	$t_{16}$          	& 1.051672 	& 0.995444	 	& 0.446226 	 		\\
 	$t_{25}$          	& 1.018130 	& 0.944031	 	& 0.412339		 	\\
  $t_{40}$          	& 1.000643 	& 0.915598	 	& 0.393605 			\\
 	$t_{41}$          	& 0.9999998	& 0.914524	 	& 0.392898 	 		\\
	$t_{100}$           & 0.986164 	& 0.890804	 	& 0.377281   		\\
\hline
%	$\NM(1.5,0.008)$      & 0.986152  & 0.887136    & 0.375221 \\
%	$\NM(1.5,0.00175)$    & 0.979775  & 0.878486    & 0.369242 \\
%	$\NM(1.5,0.000309)$   & 0.978235  & 0.876417    & 0.367833  \\
	$\NM(3,0.008)$        & 1.399511  & 1.351120    & 0.627943 \\
	$\NM(3,0.00175)$      & 1.095255  & 0.999755    & 0.429469  \\
	$\NM(3,0.000309)$     & 0.999901  & 0.898767    & 0.378687  \\
\end{tabular}
\label{tab:numeric.2}
\end{table}
%
%
%
% 
%  *=======================================================================================*
%  [\       \       \       \       \       \       \       \       \       \       \      ]
%  [ \       \       \       \       \       \       \       \       \       \       \     ]
%  [  \       \       \       \       \       \       \       \       \       \       \    ]
%  *=======================================================================================*

\subsection{The normal mixture distribution}
\label{subsec:nm}
The normal mixture distribution $\NM(\lambda,\epsilon)$, sometimes also referred to as contaminated normal distribution, is defined as
\[
	\NM(\lambda,\epsilon) \ = \ (1-\epsilon) N(0,1) + \epsilon N(0,\lambda^2), \qquad 0 \le \epsilon \le 1, \lambda \ge 1. 
\]
The resulting density is given in Table~\ref{tab:specific.1}. The normal mixture distribution is a popular model in robust statistics. It captures the notion that the majority of the data stems from the normal distribution, except for some small fraction $\epsilon$ which stems from another, usually heavier-tailed, contamination distribution. In case of the normal mixture model, this contamination distribution is the Gaussian distribution with standard deviation $\lambda$. This type of contamination model has been popularized by \citet{Tukey1960}, who also argues that $\lambda = 3$ is a sensible choice in practice. 

It is sufficient to consider the case $\lambda \ge 1$, since the parameter pair $(\lambda,\epsilon)$ yields (up to scale) the same distribution as $(1/\lambda,1-\epsilon)$. 
Now, letting $\lambda > 1$, the case where $\epsilon$ is small is the interesting one. In this case the contamination is heavy tailed, which strongly affects the behavior of our scale measures. The case $\epsilon$ close to 1 is of lesser interest: it corresponds to a normal distribution with a contamination concentrated at the origin, which affects the scale measures to a much lesser extent. 

\begin{figure}[t]
\centering
	\includegraphics[width=1.0\textwidth]{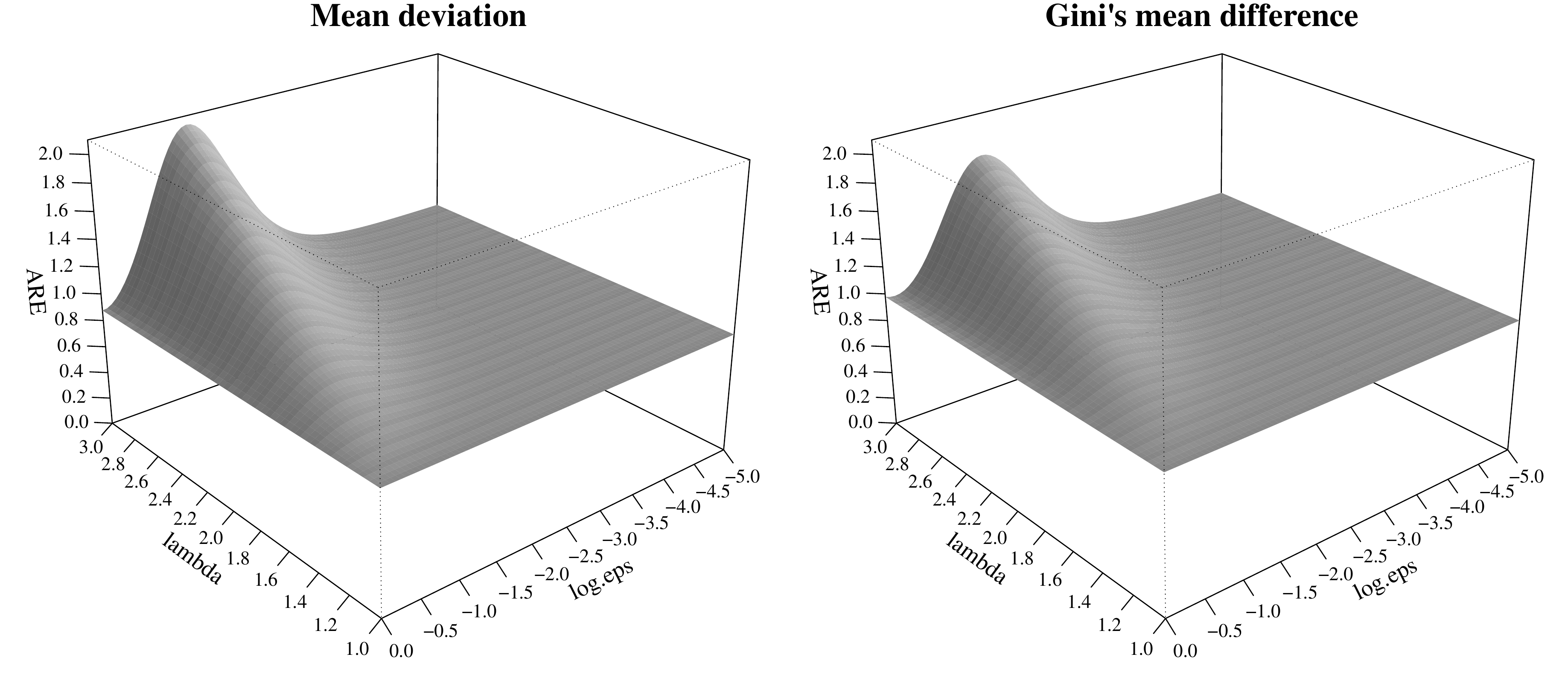}
	\includegraphics[width=0.9\textwidth]{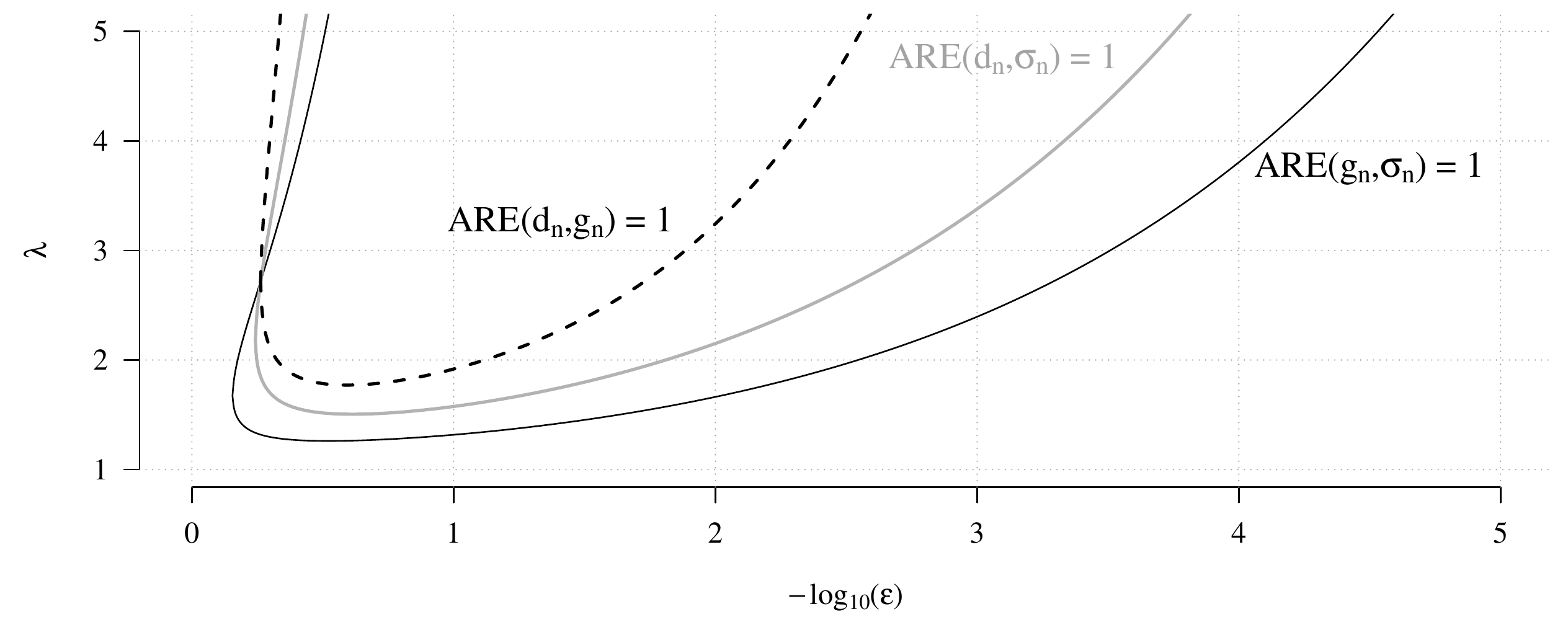}
\caption{
	Top row: Asymptotic relative efficiencies of the mean deviation (left) and Gini's mean difference (right) in the normal mixture model as a function   of $\lambda$ and $\log(\epsilon)$. 	Bottom: The curves for which values of $\lambda$ and $\epsilon$ the scale measures have the same asymptotic efficiency.
}
\label{fig:twoAREs}
\end{figure}

From the expressions for $\sigma$, $d$ and the corresponding asymptotic variances, as given in Table~\ref{tab:specific.2}, we obtain the asymptotic relative efficiency $ARE(d_n,\sigma_n)$ 
%of the mean deviation with respect to the standard deviation 
as a function of $\lambda$ and $\epsilon$. This function is plotted in Figure~\ref{fig:twoAREs} (top left). The parameter $\epsilon$ is on a log-scale since we are primarily interested in small contamination fractions. Fixing $\lambda = 3$, we find that for $\epsilon = 0.00175$, the mean deviation is as efficient as the standard deviation. \citet{Tukey1960} gives a value of $\epsilon = 0.008$. The more precise value of $0.00175$ is also in line with the simulation results of Section~\ref{sec:finite-sample}, and it supports even more so Tukey's main message: that this values is surprisingly low. 
\citet{Tukey1960} also points out that it is virtually impossible 
%for a reasonable sample size 
to distinguish a normal sample from a sample generated by a normal mixture distribution with such a low contamination fraction.

As for Gini's mean difference, the asymptotic relative efficiency $ARE(g_n,\sigma_n)$ 
%as a function of $\lambda$ and $\epsilon$ 
is depicted in the upper right plot of Figure~\ref{fig:twoAREs}. For $\lambda= 3$, Gini's mean difference is as efficient as the standard deviation for $\epsilon$ as small as $0.000309$. In the lower plot of Figure~\ref{fig:twoAREs}, equal-efficiency curves are drawn. They represent those parameter values $(\lambda,\epsilon)$, for which each two of the scale measures have equal asymptotic efficiency. So for instance, the solid black line corresponds to the contour line at height 1 of the 3D surface depicted in the top right plot. 
%
%
%  *=======================================================================================*
%  [\       \       \       \       \       \       \       \       \       \       \      ]
%  [ \       \       \       \       \       \       \       \       \       \       \     ]
%  [  \       \       \       \       \       \       \       \       \       \       \    ]
%  [   \       \       \       \       \       \       \       \       \       \       \   ]
%  [    \       \       \       \       \       \       \       \       \       \       \  ]
%  [     \       \       \       \       \       \       \       \       \       \       \ ]
%  [      \       \       \       \       \       \       \       \       \       \       \]
%  [       \       \       \       \       \       \       \       \       \       \       ]
%  *=======================================================================================*

\section{Influence functions}
\label{sec:if}
\begin{figure}[ht]
\centering
	\includegraphics[width=0.7\textwidth]{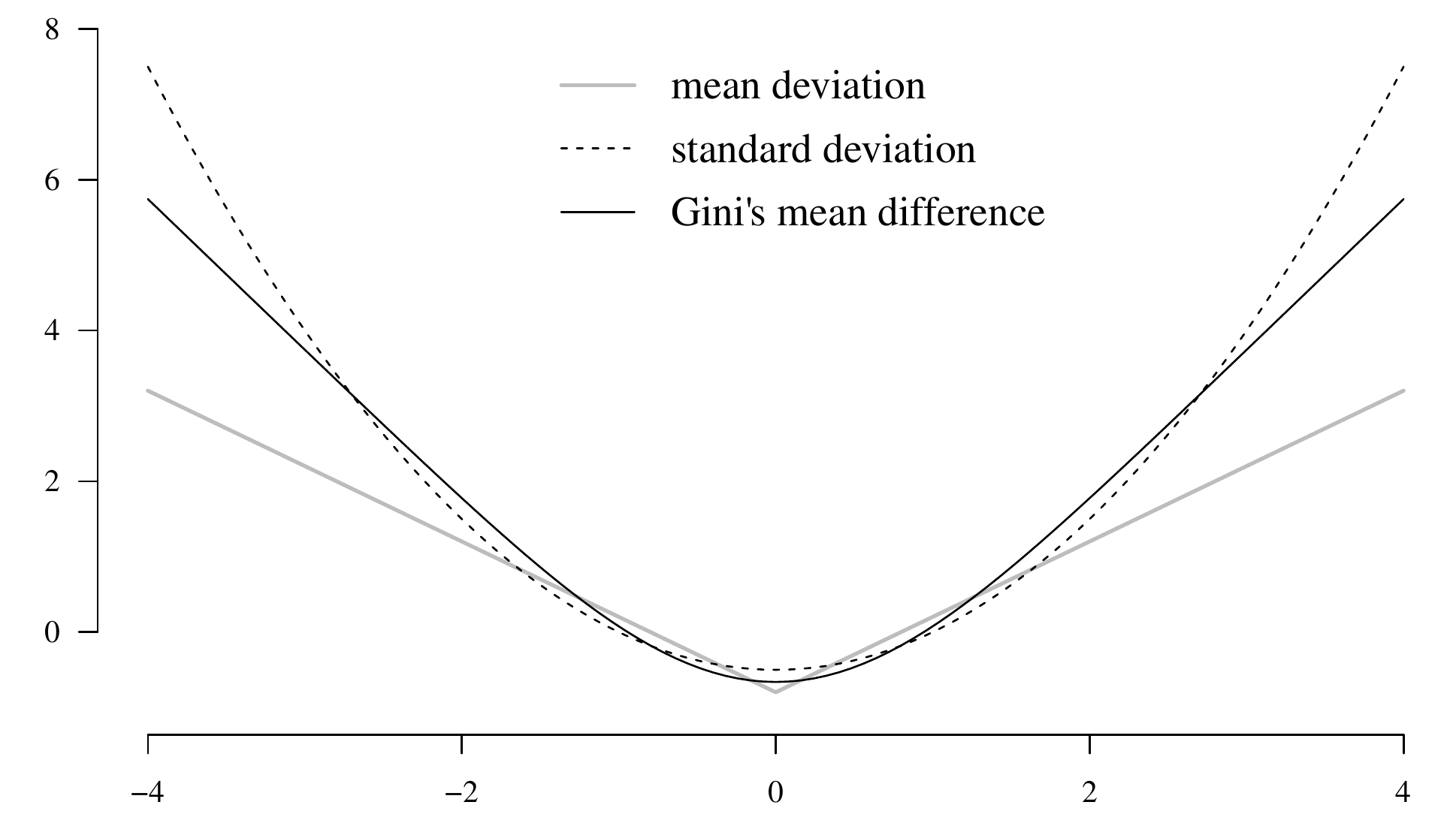}
\caption{Influence functions of the standard deviation, the mean deviation and the Gini mean difference at the standard normal distribution.}
\label{fig:if}
\end{figure}

The influence function $IF(\cdot,s,F)$ of a statistical functional $s$ at distribution $F$ is defined as
\[
	IF(x,s,F) = \lim_{\epsilon \searrow 0} \frac{1}{\epsilon} \{ s(F_{\epsilon,x}) - s(F) \},
\]
where $F_{\epsilon,x} = (1-\epsilon)F + \epsilon \Delta_x$, $0 \le \epsilon \le 1$, $x \in \R$, and $\Delta_x$ denotes Dirac's delta, i.e., the probability measure that puts unit mass in $x$. The influence function describes the impact of an infinitesimal contamination at point $x$ on the functional $s$ if the latter is evaluated at distribution $F$. For further reading see, e.g., \citet{Huber2009} or \citet{Hampel1986}.
The influence functions  of the three scale measures are
\[
	IF(x,\sigma(\cdot);F) \ = \ (2\sigma(F))^{-1}\{ (E(X) - x)^2 - \sigma^2(F) \}, 
\]\[
	IF(x,d(\cdot); F) \ = \ |x-m(F)| - d(F), 
\]\[
	IF(x,g(\cdot); F) \ = \ 2 \left\{ x [ F(x) + F(x-) - 1 ] + E[X \ind{X\ge x}] -  E[X \ind{X \le x}]  - g(F) \right\}
\]
The derivations are straightforward, and the results are stated without proof. 
For the formula of $d(\cdot)$ to hold, $F$ has to fulfill certain regularity conditions in the vicinity of its median $m(F)$. Specifically,
$(m(F_{\epsilon,x})-m(F)) = O(\epsilon)$ as $\epsilon \to 0$ for all $x \in \R$ and $F(m(F_{\epsilon,x})) \to 1/2$ are a set of sufficient conditions. They are fulfilled, e.g., if $F$ possesses a positive Lebesgue density in a neighborhood of $m(F)$. 
For the standard normal distribution, above expressions reduce to
\[
	IF(x,\sigma(\cdot);N(0,1)) \ = \ (x^2-1)/2, 
\]\[
	IF(x,d(\cdot); N(0,1)) \ = \ |x| - \sqrt{2/\pi}, 
\]\[
	IF(x,g(\cdot); N(0,1)) \ = \ 4\phi(x) +  2 x \{ 2\Phi(x) - 1 \} - 4/\sqrt{\pi}, 
\]
where $\phi$ and $\Phi$ denote the density and the cdf of the standard normal distribution, respectively.
These curves are depicted in Figure \ref{fig:if}. They confirm the general impression mediated by Table \ref{tab:numeric.2}, that Gini's mean difference is in-between the standard and the mean deviation, and support our claim that it combines the advantages of the other two: its influence function grows linearly for large $|x|$, but it is smooth at the origin.   

%
%
%  *=======================================================================================*
%  [\       \       \       \       \       \       \       \       \       \       \      ]
%  [ \       \       \       \       \       \       \       \       \       \       \     ]
%  [  \       \       \       \       \       \       \       \       \       \       \    ]
%  [   \       \       \       \       \       \       \       \       \       \       \   ]
%  [    \       \       \       \       \       \       \       \       \       \       \  ]
%  [     \       \       \       \       \       \       \       \       \       \       \ ]
%  [      \       \       \       \       \       \       \       \       \       \       \]
%  [       \       \       \       \       \       \       \       \       \       \       ]
%  *=======================================================================================*

\section{Finite sample efficiencies}
\label{sec:finite-sample}

In a small simulation study we want to check if the asymptotic efficiencies computed in Section \ref{sec:asymptotic}
are useful approximations for the actual efficiencies in finite samples. For this purpose we consider the following nine distributions: the standard normal $N(0,1)$, the standard Laplace $L(0,1)$ (with parameters $\mu =0$ and $\alpha = 1$, cf.~Table~\ref{tab:specific.1}), the uniform distribution $U(0,1)$ on the unit interval, the $t_\nu$ distribution with $\nu= 5, 16, 41$ and the normal mixture with the parameter choices as in Tables~\ref{tab:numeric.1} and \ref{tab:numeric.2}. The choice $\nu = 5$ serves as a heavy-tailed example, whereas for $\nu = 16$ and $\nu = 41$ we have witnessed at Table~\ref{tab:numeric.2} that the mean deviation and the Gini mean difference, respectively, are asymptotically equally efficient as the standard deviation.

For each distribution and each of the sample sizes $n = 5, 8, 10, 50, 500$, we generate 100,000 samples and compute from each sample the three scale measures $\sigma_n$, $d_n$ and $g_n$.  The results for $N(0,1)$, $L(0,1)$ and $U(0,1)$ are summarized in Table~\ref{tab:finite-sample.1}, for the $t_\nu$ distributions in Table~\ref{tab:finite-sample.2} and for the normal mixture distributions in Table~\ref{tab:finite-sample.3}. 
For each estimate, population distribution and sample size, the following numbers are reported: the sample variance of the 100,000 estimates multiplied by the respective value of $n$ (the ``$n$-standardized variance'' which approaches the asymptotic variance given in Table~\ref{tab:numeric.1} as $n$ increases), the squared bias relative to the variance, and, for $g_n$ and $d_n$, the relative efficiencies with respect to the standard deviation. For the relative efficiency computation, it is important to note that the standardizing, cf.~(\ref{eq:are}), is done not by the true asymptotic values, but by the empirical finite-sample value, i.e. the sample mean of the 100,000 estimates. Also note that variances, not mean squared errors, are reported.   
For Gini's mean difference, the simulated variances are also compared to the true finite-sample variances, cf.~(\ref{eq:var.g_n}).

%From Tables \ref{tab:finite-sample.1} to \ref{tab:finite-sample.3}, 
We observe the following:
For large and moderate sample sizes ($n = 50, 500$), the simulated values are close to the asymptotic ones from Tables \ref{tab:numeric.1} and \ref{tab:numeric.2}, and we conclude that the asymptotic efficiency generally provides a useful indication for the actual efficiency.
In small samples, the simulated relative efficiencies may substantially differ from the asymptotic values, but the ranking of the three estimators stays the same. %However, the small-sample efficiency comparison should be treated with caution since there is a noticeable bias of the standard deviation. % eigentlich blÖdsinn, nochmal besser schreiben
Furthermore, the simulations confirm the unbiasedness of Gini's mean difference and the formula (\ref{eq:var.g_n}), due to \citet{lomnicki:1952}, for its finite-sample variance.

Finally, we also include the mean deviation with factor $1/n$ instead of $1/(n-1)$ in the study, denoted by $d_n^*$ in the tables. Since $d_n$ and $d_n^*$ differ only by multiplicative factor, the efficiencies are the same, and we only report the (squared) bias (relative to the variance). We find that $d_n^*$ is heavily biased for small samples for all distributions considered, whereas $d_n$ has in all situations a smaller bias than $\sigma_n$. 
Somewhat unexpected is the increase of the squared bias relative to the variance of $d_n$ from $n = 5$ to $n = 8$. The reason may lie in the different behavior of the sample median for odd and even numbers of observations. 

The simulations were done in R \citep{R}, using an implementation for Gini's mean difference by A.~Azzalini.\footnote{https://stat.ethz.ch/pipermail/r-help/2003-April/032820.html}

\begin{table}
\caption{Simulated variances, biases and relative efficiencies of $\sigma_n$, $g_n$, $d_n$ at $N(0,1)$, $L(0,1)$ and $U(0,1)$ for several sample sizes, $d_n^*$: mean deviation with $1/n$ scaling.}
\small
\label{tab:finite-sample.1}
%\centering
\begin{tabular}{c|c|r|r|r|r|r}
\hline%\hline
 estimator 	&  													& $n=5$ 		& $n=8$ 	 & $n=10$ 	& $n=50$ 	 & $n=500$	\\
\hline\hline 

\multicolumn{2}{c}{ \parbox[][5.0ex][c]{0.15\textwidth}{$N(0,1)$}} &\multicolumn{5}{c}{} \\ \hline
$\sigma_n$ 	&$n \cdot$variance					  &0.577	&0.548 &0.541 &0.507 &0.505 			\\
						&bias$^2/$variance						&0.031	&0.019 &0.014 &0.003 &$<$ 0.001		\\
\hline
$g_n$ 			&$n \cdot$variance (empirical)&0.850	&0.767 &0.743 &0.666 &0.655 			\\
						&$n \cdot$variance (true) 		&0.852  &0.766 &0.740 &0.667 &0.653				\\
						&bias$^2/$variance 						&3.4e-08	&4.7e-07 &7.8e-06 &1.0e-05 &4.7e-06				\\
	&rel.\ efficiency\ wrt $\sigma_n$				&0.986	&0.982 &0.980 &0.979 &0.978				\\
\hline
$d_n$			&$n\cdot$variance     				  &0.482  &0.454 &0.427 &0.374 &0.365 			\\
					&bias$^2/$variance 							&0.009	&0.020 &0.012 &0.001 &$<$ 0.001		\\
			&rel.\ efficiency\ wrt $\sigma_n$		&0.938	&0.902 &0.894 &0.880 &0.876				\\
\hline $d_n^*$				&bias$^2/$variance    			&0.296	&0.118 &0.101 &0.021 &0.002				\\
\hline\hline 
\multicolumn{2}{c}{\parbox[][5.0ex][c]{0.15\textwidth}{$L(0,1)$}}&\multicolumn{5}{c}{}\\
\hline
$\sigma_n$ 	&$n\cdot$variance					  	&1.946  &2.076 &2.134 &2.387 &2.495					\\
						&bias$^2/$variance						&0.055	&0.034 &0.027 &0.006 &$<$0.001			\\
\hline
$g_n$ 			&$n\cdot$variance (empirical)	&2.629  &2.514 &2.456 &2.359 &2.345 				\\
						&$n\cdot$variance (true) 			&2.625  &2.500 &2.463 &2.357 &2.336					\\
						&bias$^2/$variance 						&2.8e-06&8.4e-09 &8.4e-08 &1.3e-05 &8.3e-10	\\
	&rel.\ efficiency\ wrt $\sigma_n$				&1.037	&1.071 &1.088 &1.167 &1.201					\\
\hline
$d_n$				&$n\cdot$variance     				&1.343  &1.232 &1.169 &1.041 &1.005 				\\
						&bias$^2/$variance 						&0.025	&0.028 &0.021 &0.005 &$<$0.001			\\
	&rel.\ efficiency\ wrt $\sigma_n$				&1.061	&1.101 &1.123 &1.206 &1.245					\\
\hline $d_n^*$				&bias$^2/$variance  &0.106	&0.040 &0.031 &0.006 &0.001					\\ 
\hline\hline 
\multicolumn{2}{c}{\parbox[][5.0ex][c]{0.15\textwidth}{$U(0,1)$}}&\multicolumn{5}{c}{}\\
\hline
$\sigma_n$ 	&$n\cdot$variance					  	&0.031 & 0.025 & 0.022 & 0.018 & 0.017 			\\
						&bias$^2/$variance						&0.021 & 0.010 & 0.007 & 0.001 & $<$ 0.001  \\
\hline
$g_n$ 			&$n\cdot$variance (empirical) &0.045 & 0.035 & 0.032 & 0.024 & 0.023 			\\
						&$n\cdot$variance (true) 	  	&0.044 & 0.035 & 0.032 & 0.024 & 0.022 			\\
						&bias$^2/$variance 					&1.9e-05&6.2e-07 &9.4e-07 &3.0e-05 &5.1e-08		\\
				&rel.\ efficiency\ wrt $\sigma_n$	&0.985 & 0.967 & 0.962 & 0.985 & 0.998      \\
\hline
$d_n$			&$n\cdot$variance               &0.030 & 0.028 & 0.026 & 0.022 & 0.021			\\
						&bias$^2/$variance 					&6.1e-06	&4.7e-03 &2.3e-03 &6.8e-05 &1.7e-05	\\
				&rel.\ efficiency\ wrt $\sigma_n$&0.829 & 0.694 & 0.672 & 0.614 & 0.603       \\
\hline $d_n^*$				&bias$^2/$variance           &0.657 & 0.285 & 0.236 & 0.059 & 0.006			\\
\hline%\hline
\end{tabular}
\end{table}

\begin{table}
\caption{Simulated variances, biases and relative efficiencies of $\sigma_n$, $g_n$, $d_n$ at $t_\nu$ distributions for several sample sizes and values of $\nu$; $d_n^*$: mean deviation with $1/n$ scaling.}
\label{tab:finite-sample.2}
\small
\begin{tabular}{c|c|r|r|r|r|r}
\hline
 estimator 	&  													& $n=5$		  & $n=8$		 & $n=10$		& $n=50$ 	 & $n=500$\\
\hline\hline
\multicolumn{2}{c}{\parbox[][5.0ex][c]{0.15\textwidth}{$t_5$}}&\multicolumn{5}{c}{} \\
\hline
$\sigma_n$ 	&$n\cdot$variance									&1.584  &1.686 &1.762 &2.313 &2.880 				\\
						&bias$^2/$variance								&0.050	&0.034 &0.028 &0.007 &0.001					\\\hline
$g_n$ 			&$n\cdot$variance (empirical)   	&2.050  &1.942 &1.890 &1.805 &1.790					\\
						&$n\cdot$variance (true) 					&2.047  &1.935 &1.901 &1.806 &1.787					\\
						&bias$^2/$variance 								&4.0e-06&1.3e-05 &2.5e-05 &5.1e-06 &1.4e-05					\\
						&rel.\ efficiency\ wrt $\sigma_n$	&1.073	&1.150 &1.185 &1.499 &1.811					\\\hline
$d_n$			&$n\cdot$variance     				   		&1.036  &0.949 &0.901 &0.791 &0.760					\\
						&bias$^2/$variance 								&0.014	&0.018 &0.014 &0.003 &$<$ 0.001 				\\
						&rel.\ efficiency\ wrt $\sigma_n$	&1.105	&1.208 &1.282 &1.673 &1.977					\\\hline
$d_n^*$				&bias$^2/$variance    					&0.160	&0.066 &0.053 &0.011 &0.001					\\
\hline\hline 
\multicolumn{2}{c}{\parbox[][5.0ex][c]{0.15\textwidth}{$t_{16}$}}&\multicolumn{5}{c}{}\\
\hline
$\sigma_n$ 	&$n\cdot$variance					  			&0.745  &0.722 &0.722 &0.710 &0.705					\\
						&bias$^2/$variance								&0.034	&0.021 &0.015 &0.003 &$<$ 0.001 	 \\
\hline
$g_n$ 			&$n\cdot$variance (empirical)   	&1.064  &0.977 &0.949 &0.862 &0.850					\\
						&$n\cdot$variance (true) 					&1.065  &0.972 &0.945 &0.866 &0.850					\\
						&bias$^2/$variance 								&1.4e-07&5.0e-06 &7.9e-07 &7.6e-06 &2.4e-05	\\
						&rel.\ efficiency\ wrt $\sigma_n$	&0.999	&1.009 &1.016 &1.043 &1.050				\\
\hline
$d_n$			&$n\cdot$variance     				   		&0.588  &0.547 &0.517 &0.454 &0.445 				\\
						&bias$^2/$variance 								&0.012	&0.018 &0.012 &0.002 &$<$ 0.001					\\
						&rel.\ efficiency\ wrt $\sigma_n$	&0.972	&0.956 &0.963 &0.989 &0.991					\\
\hline $d_n^*$				&bias$^2/$variance    					&0.259	&0.106 &0.085 &0.017 &0.002					\\
\hline\hline 
\multicolumn{2}{c}{\parbox[][5.0ex][c]{0.15\textwidth}{$t_{41}$}}&\multicolumn{5}{c}{}\\
\hline
$\sigma_n$ 	&$n\cdot$variance					  			&0.640  &0.611 &0.605 &0.574 &0.575					\\
						&bias$^2/$variance								&0.032	&0.020 &0.014 &0.003 &$<$ 0.001			\\
\hline
$g_n$ 			&$n\cdot$variance (empirical)   	&0.925  &0.835 &0.817 &0.740 &0.720 				\\
						&$n\cdot$variance (true) 					&0.925  &0.837 &0.811 &0.736 &0.720					\\
						&bias$^2/$variance 								&1.1e-05	&3.6e-06 &1.5e-06 &9.5e-08 &7.1e-07	\\
						&rel.\ efficiency\ wrt $\sigma_n$	&0.990	&0.992 &0.991 &0.999 &1.001				\\
\hline
$d_n$			&$n\cdot$variance     				   		&0.519  &0.482 &0.462 &0.399 &0.390				\\
						&bias$^2/$variance 								&0.010	&0.018 &0.013 &0.002 &$<$ 0.001		\\
						&rel.\ efficiency\ wrt $\sigma_n$	&0.950	&0.918 &0.921 &0.916 &0.919					\\
\hline $d_n^*$				&bias$^2/$variance    					&0.276	&0.113 &0.094 &0.019 &0.002					\\
\hline
\end{tabular}
\end{table}

\begin{table}
\caption{Simulated variances, biases and relative efficiencies of $\sigma_n$, $g_n$, $d_n$ at normal mixture distributions for $\lambda=3$ and $\epsilon = 0.008, 0.00175, 0.000309$; $d_n^*$: mean deviation with $1/n$ scaling.}
\label{tab:finite-sample.3}
\small
\begin{tabular}{c|c|r|r|r|r|r}
\hline
 estimator 	&  													& $n=5$		  & $n=8$		 & $n=10$		& $n=50$ 	 & $n=500$\\
\hline\hline
\multicolumn{2}{c}{\parbox[][5.0ex][c]{0.15\textwidth}{$\NM(3,0.008)$}}&\multicolumn{5}{c}{} \\
\hline
$\sigma_n$ 	&$n\cdot$variance					  	&0.710  &0.698 &0.711 &0.815 &0.875					\\
						&bias$^2/$variance						&0.034	&0.024 &0.018 &0.004 &0.001					\\
\hline
$g_n$ 			&$n\cdot$variance (empirical) &0.997  &0.910 &0.876 &0.804 &0.790					\\
						&$n\cdot$variance (true) 			&0.996  &0.908 &0.882 &0.808 &0.793					\\
						&bias$^2/$variance 						&4.6e-06&2.1e-10 &1.6e-05 &3.4e-06 &8.4e-07	\\
			&rel.\ efficiency\ wrt $\sigma_n$		&1.023	&1.060 &1.083 &1.257 &1.385					\\
\hline
$d_n$				&$n\cdot$variance     				&0.540  &0.507 &0.480 &0.423 &0.405 				\\
						&bias$^2/$variance 						&0.010	&0.016 &0.013 &0.002 &$<$ 0.001					\\
		&rel.\ efficiency\ wrt $\sigma_n$			&1.000	&1.016 &1.039 &1.204 &1.332					\\
\hline $d_n^*$			&bias$^2/$variance    &0.264	&0.112 &0.087 &0.020 &0.002					\\
\hline\hline
\multicolumn{2}{c}{\parbox[][5.0ex][c]{0.15\textwidth}{$\NM(3,0.00175)$}}&\multicolumn{5}{c}{} \\
\hline 
$\sigma_n$ 	&$n\cdot$variance					  	&0.617  &0.587 &0.576 &0.573 &0.590					\\
						&bias$^2/$variance						&0.032	&0.019 &0.017 &0.003 &$<$ 0.001			\\
\hline
$g_n$ 			&$n\cdot$variance (empirical) &0.889  &0.791 &0.764 &0.704 &0.675					\\
						&$n\cdot$variance (true) 			&0.883  &0.797 &0.771 &0.698 &0.684					\\
						&bias$^2/$variance 						&1.6e-07&3.0e-07 &4.8e-08 &1.8e-05 &1.0e-05	\\
				&rel.\ efficiency\ wrt $\sigma_n$	&0.995	&1.002 &1.009 &1.056 &1.092					\\
\hline
$d_n$				&$n\cdot$variance     				&0.500  &0.462 &0.441 &0.385 &0.370					\\
						&bias$^2/$variance 						&0.011  &0.017 &0.013 &0.002 &3.9e-05				\\
				&rel.\ efficiency\ wrt $\sigma_n$	&0.951	&0.931 &0.931 &0.971 &0.992					\\
\hline $d_n^*$			&bias$^2/$variance    				&0.283	&0.115 &0.100 &0.022 &0.003					\\
\hline\hline
\multicolumn{2}{c}{\parbox[][5.0ex][c]{0.15\textwidth}{$\NM(3,0.000309)$}}&\multicolumn{5}{c}{} \\
\hline 
	$\sigma_n$ 	&$n\cdot$variance					  &0.584  &0.558 &0.543 &0.517 &0.515			\\
						&bias$^2/$variance						&0.031	&0.017 &0.014 &0.003 &$<$ 0.001 \\
\hline
$g_n$ 			&$n\cdot$variance (empirical) &0.853  &0.775 &0.744 &0.667 &0.655			\\
						&$n\cdot$variance (true) 			&0.857  &0.771 &0.746 &0.673 &0.658			\\
						&bias$^2/$variance 						&1.3e-05	&4.8e-06 &5.1e-07 &1.3e-06 &8.3e-06	\\
			  &	rel.\ efficiency\ wrt $\sigma_n$	&0.986	&0.986 &0.985 &0.993 &0.999			\\
\hline
$d_n$				&$n\cdot$variance     				&0.484  &0.452 &0.434 &0.375 &0.365 		\\
						&bias$^2/$variance 						&0.009  &0.018 &0.012 &0.002 &$<$ 0.001  \\
				&rel.\ efficiency\ wrt $\sigma_n$	&0.941  &0.900 &0.903 &0.899 &0.903 		\\
\hline $d_n^*$			&bias$^2/$variance    				&0.291	&0.122 &0.096 &0.021 &0.002			\\
\hline
\end{tabular}
\end{table}

%
%
%  *=======================================================================================*
%  [\       \       \       \       \       \       \       \       \       \       \      ]
%  [ \       \       \       \       \       \       \       \       \       \       \     ]
%  [  \       \       \       \       \       \       \       \       \       \       \    ]
%  [   \       \       \       \       \       \       \       \       \       \       \   ]
%  [    \       \       \       \       \       \       \       \       \       \       \  ]
%  [     \       \       \       \       \       \       \       \       \       \       \ ]
%  [      \       \       \       \       \       \       \       \       \       \       \]
%  [       \       \       \       \       \       \       \       \       \       \       ]
%  *=======================================================================================*

\section{Summary and conclusion}
\label{sec:conclusion}
Neither the standard deviation nor the mean deviation is a robust estimator. However, several authors have argued that, when comparing the standard deviation with the mean deviation, the (relatively) better robustness of the latter
%, its reduced vulnerability to outliers, heavy-tailed distributions or simply any departures from normality, 
is a crucial advantage, which outweighs its disadvantages, and that the mean deviation is hence to be preferred out of the two. We share this view. However, we recommend to use Gini's mean difference instead of the mean deviation. While it has qualitatively the same robustness and the same efficiency under long-tailed distributions as the mean deviation% 
%and requires also only second moments for asymptotic normality
, it lacks its main disadvantage as compared the standard deviation: the lower efficiency at strict normality. For near-normal distributions --- and also for very light-tailed distribution, as the results for the uniform distribution suggest ---,
Gini's mean difference and the standard deviation are for all practical purposes equally efficient. For instance, at the normal and all $t_\nu$ distributions with $\nu \ge 23$, the (properly standardized) asymptotic variances of $g_n$ and $\sigma_n$ are within a three percent margin of each other. At heavy-tailed distributions, Gini's mean difference is, along with the mean deviation, substantially more efficient than the standard deviation. However, it must also be noted that heavy tails are a bad case scenario for all three scale measure, and that, for heavy-tailed distributions, much more efficient estimators available, e.g., $M$-estimators or trimmed or quantile-based scale estimators.

Gini's mean difference has further advantages, which particularly concern its finite-sample performance: it is unbiased, and the finite-sample variance is known. It either can be computed exactly or, if no specific model is assumed, estimated, which allows for instance better approximative confidence intervals. Neither of that is true for the standard deviation or the mean deviation, and one can consequently argue that Gini's mean difference is a superior scale estimator even under normality.
Scale measures may serve different purposes: some, such as the interquartile range, are primarily used for descriptive reasons. In this respect, the mean deviation and the standard deviation may be preferred (the latter due to its widespread use), but as an estimator, i.e, for inferring about an unknown population scale, Gini's mean difference has, in our opinion, clear advantages. Although there are quite a few articles that advocate the use of alternative scale estimators \citep[specifically for Gini's mean difference, cf.\ e.g.][]{Yitzhaki2003},  we are aware that the standard deviation is so widespread and common 
%as \emph{the} standard measure of dispersion 
that any other scale estimator taking its role seems as unlikely as a change from the decimal to another numeral system. 

%
%
%  *=======================================================================================*
%  [\       \       \       \       \       \       \       \       \       \       \      ]
%  [ \       \       \       \       \       \       \       \       \       \       \     ]
%  [  \       \       \       \       \       \       \       \       \       \       \    ]
%  [   \       \       \       \       \       \       \       \       \       \       \   ]
%  [    \       \       \       \       \       \       \       \       \       \       \  ]
%  [     \       \       \       \       \       \       \       \       \       \       \ ]
%  [      \       \       \       \       \       \       \       \       \       \       \]
%  [       \       \       \       \       \       \       \       \       \       \       ]
%  *=======================================================================================*

\section*{Acknowledgment}
We are indebted to Herold Dehling for introducing us to the theory of $U$-statistics, to Roland Fried for introducing us to robust statistics, and to Alexander D\"urre, who has demonstrated the benefit of complex analysis for solving statistical problems. Both authors were supported in part by the Collaborative Research Centre 823 \emph{Statistical modelling of nonlinear dynamic processes}.

%
%
%  *=======================================================================================*
%  [\       \       \       \       \       \       \       \       \       \       \      ]
%  [ \       \       \       \       \       \       \       \       \       \       \     ]
%  [  \       \       \       \       \       \       \       \       \       \       \    ]
%  [   \       \       \       \       \       \       \       \       \       \       \   ]
%  [    \       \       \       \       \       \       \       \       \       \       \  ]
%  [     \       \       \       \       \       \       \       \       \       \       \ ]
%  [      \       \       \       \       \       \       \       \       \       \       \]
%  [       \       \       \       \       \       \       \       \       \       \       ]
%  *=======================================================================================*

\appendix
\section{Integrals for the normal distribution}
\label{sec:app:norm}

\begin{figure}[ht]
\centering
	\includegraphics[width=0.4\textwidth]{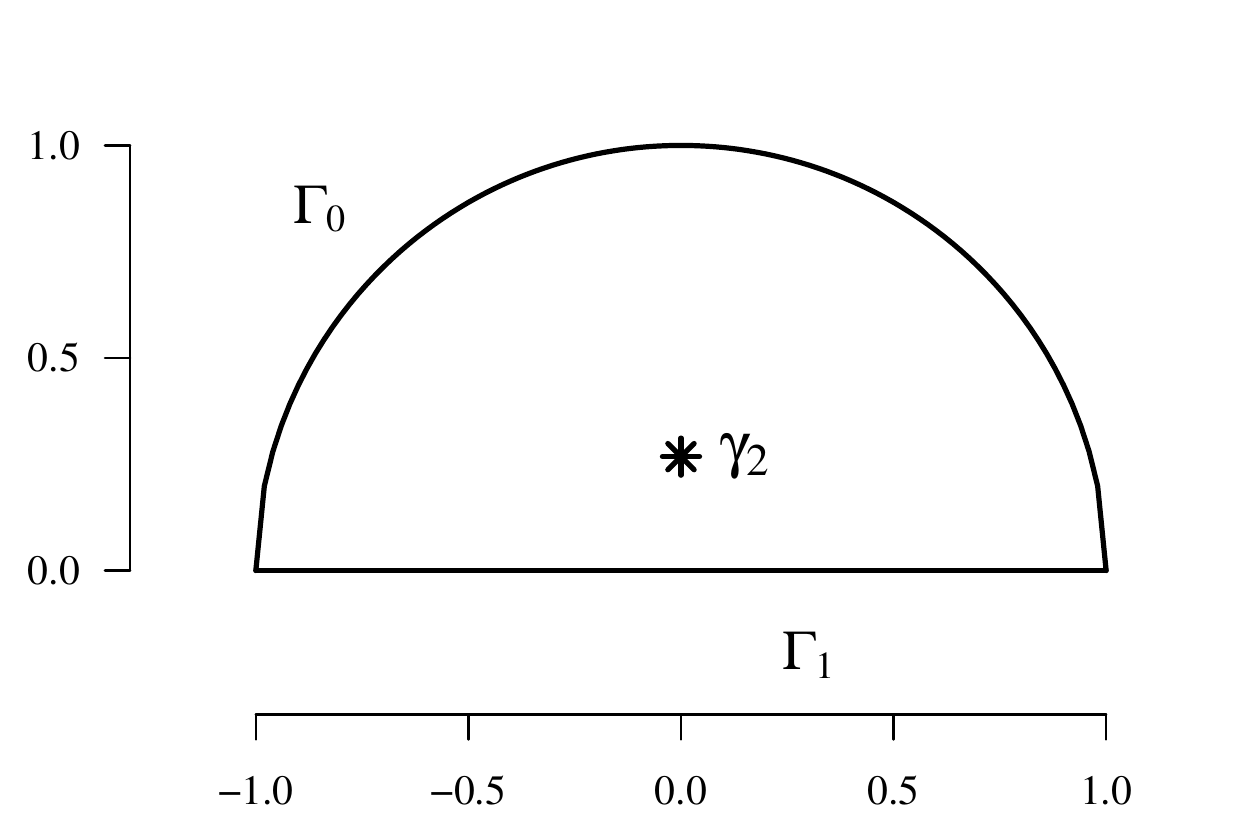}
\caption{Residue theorem: the line integral over $h$ along the closed curve $\Gamma=\Gamma_0\, \cup\, \Gamma_1$ is determined by the residue of $h$ in $\gamma_2$.}
\label{fig:complex}
\end{figure}

When evaluating the integral $J$, cf.~(\ref{eq:J}), for the standard normal distribution, one encounters the integral 
\[
	I_1 = \int_{-\infty}^{\infty} x^2 \phi(x) \Phi(x)^2 dx, 
\]
where $\phi$ and $\Phi$ denote the density and the cdf of the standard normal distribution, respectively. \citet{nair:1936} gives the value 
$I_1 = 1/3 + 1/(2\pi\sqrt{3})$, resulting in $J = \sqrt{3}/(2 \pi) - 1/6$, but does not provide a proof. The author refers to the derivation of a similar integral \citep[integral 8 in Table I,][p.~433]{nair:1936}, where we find the result as well as the derivation doubtful, and to an article by \citet{Hojo:1931}, which gives numerical values for several integrals, but does not contain an explanation for the value of $I_1$ either. We therefor include a proof here. 
Writing  $\Phi(x)$ as the integral of its density and changing the order of the integrals in thus obtained three-dimensional integral yields
\[
	I_1 = 
	(2\pi)^{-3/2}
		\int_{y=-\infty}^0 \int_{z=\infty}^0 \int_{x=-\infty}^\infty  
		x^2 e^{x^2/2} e^{(y+x)^2/2} e^{(z+x)^2/2} 
	\, d x\, d z\, d y.
\]
Solving the inner integral, we obtain
\[	
	I_1 = 
	(18\pi\sqrt{3})^{-1} 
	\int_{y=0}^\infty \int_{z=0}^\infty  [ (y+z)^2 + 3 ] 
	\exp\left\{ - \frac{1}{3} \left[  y^2 + z^2 - y z \right]\right\}
	\, d z\, d y.
\]
Introducing polar coordinates $\alpha, r$ such that $y = r \cos \alpha$, $z = r \sin \alpha$, and solving the integral with respect to $r$, we arrive at
\[
	I_1 = 	
	\frac{1}{4\pi\sqrt{3}} 
	\int_{\alpha=0}^\pi \frac{4+\sin \alpha}{(2-\sin \alpha)^2} \, d \alpha.
\]
This remaining integral may be solved by means of the residue theorem \citep[e.g.][p.~149]{Ahlfors1966}. Substituting $\gamma = e^{i \alpha}$ and using $\sin \alpha= (e^{i \alpha} - e^{- i \alpha})/(2i)$, we transform $I_1$ into the following line integral in the complex plane, 
\be \label{eq:complex}
	I_1 = 
	\frac{1}{4\pi\sqrt{3}} 
	\int_{\Gamma_0} \frac{\gamma^2 + 8 i \gamma  -1 }{(\gamma^2- 4 i \gamma  -1 )^2} \, d \gamma,
\ee
where $\Gamma_0$ is the upper unit half circle in the complex plane, cp.~Figure \ref{fig:complex}. 
Let us call $h$ the integrand in (\ref{eq:complex}), its poles (both of order two) are $\gamma_{1/2} = (2\pm\sqrt{3})i$, so that $\gamma_2$ lies within the closed upper half unit circle $\Gamma$. The residue of $h$ in $\gamma_2$ is $-\sqrt{3} i /2$. Integrating $h$ along $\Gamma_1$, i.e.~the real line from -1 to 1, cf.~Figure \ref{fig:complex}, and applying the residue theorem to the closed line integral along $\Gamma$ completes the derivation. 

%
%
%  *=======================================================================================*
%  [\       \       \       \       \       \       \       \       \       \       \      ]
%  [ \       \       \       \       \       \       \       \       \       \       \     ]
%  [  \       \       \       \       \       \       \       \       \       \       \    ]
%  [   \       \       \       \       \       \       \       \       \       \       \   ]
%  [    \       \       \       \       \       \       \       \       \       \       \  ]
%  [     \       \       \       \       \       \       \       \       \       \       \ ]
%  [      \       \       \       \       \       \       \       \       \       \       \]
%  [       \       \       \       \       \       \       \       \       \       \       ]
%  *=======================================================================================*

\section{Integrals for the normal mixture distribution}

Evaluating the integral $J$ for the normal mixture distribution, we arrive after some lengthy but straightforward calculations at
\[
J \ = \ \Big[ \epsilon^3 \lambda^2 + (1-\epsilon)^3 \Big] \Big[ 2 A(1) + C(1) + E(1) \Big] \ - \ (\epsilon \lambda^2 + 1 - \epsilon)B 
\]\[
	\qquad + \ \epsilon^2(1-\epsilon) 
		\Big[
				2 (2+ \lambda^2) A(1/\lambda) + C(\lambda) + 2\lambda^2 D(1/\lambda)+ \lambda(2+ \lambda^2)E(1/\lambda) 
		\Big]
\]\[
	\qquad + \ \epsilon (1-\epsilon)^2
		\Big[
				2 (2\lambda^2+ 1 ) A(\lambda) + \lambda^2C(1/\lambda) + 2 D(\lambda)+  (\lambda^{-1} + 2 \lambda )E(\lambda) 
		\Big],
\]
where 
\[
 	A(\lambda) = \int_{-\infty}^\infty x \phi^2(x) \Phi(x/\lambda) dx \ = \ \frac{1}{4 \pi \sqrt{ 1 + 2 \lambda^2 }},
 	\qquad
 	B  = \int_{-\infty}^\infty x^2 \phi(x) \Phi(x) dx \ = \ \frac{1}{2},
\]
{\small \[
	C(\lambda) = \int_{-\infty}^\infty x^2 \phi(x) \Phi^2(x/\lambda) dx 
	\ = \ \frac{1}{4} 
				+ \frac{\lambda}{\pi (1 + \lambda^2)\sqrt{ 2 + \lambda^2 }} 
				+ \frac{1}{2 \pi} \arctan\left(\frac{1}{\lambda\sqrt{2+\lambda^2}}\right),
\] \[
	D(\lambda) = \int_{-\infty}^\infty x^2 \phi(x) \Phi(x) \Phi(x/\lambda) dx 
	\ = \ \frac{1}{4} 
				+ \frac{3\lambda^2 + 1}{4 \pi  (1 + \lambda^2)\sqrt{ 2 \lambda^2 + 1 }} 
				+ \frac{1}{2 \pi} \arctan\left(\frac{1}{\sqrt{2\lambda^2+ 1}}\right),
\]} % small endet
\[
	E(\lambda) = \int_{-\infty}^\infty  \phi^2(x) \phi(x/\lambda) dx \ = \ \frac{1}{2 \pi \sqrt{ 1 + 2 \lambda^2 }},
\]
for all $\lambda > 0$. As before, $\phi$ and $\Phi$ denote the density and the cdf of standard normal distribution.
The tricky integrals are $C(\lambda)$ and $D(\lambda)$, which, for $\lambda = 1$, both reduce to the integral $I_1$ above. They can be solved by similar means as $I_1$. Proceeding as in Appendix \ref{sec:app:norm}, solving the respective two inner integrals yields
\[
	C(\lambda) = \frac{\lambda^3}{2 \pi \sqrt{2+\lambda^2}} \int_0^{\pi/2} \frac{ 3 + \lambda^2 + \sin(2\alpha)}{ \{1 + \lambda^2 - \sin(2 \alpha)\}^2} d\alpha,
\]\[
	D(\lambda) = \frac{1}{2 \pi\sqrt{1+2\lambda^2}}  \int_0^{\pi/2} 
	\frac{ 2 + \lambda^2 (2+ \sin(2\alpha)) + (3\lambda^4 - \lambda^2 - 2) \sin^2(\alpha)}
	      { \{ 2 - \sin(2\alpha) + (\lambda^2-1) \sin^2(\alpha) \}^2 } d \alpha.
\]
These integrals are again solved by the residue theorem, for which we used the software Mathematica \citep{Mathematica}.

%
%
%  *=======================================================================================*
%  [\       \       \       \       \       \       \       \       \       \       \      ]
%  [ \       \       \       \       \       \       \       \       \       \       \     ]
%  [  \       \       \       \       \       \       \       \       \       \       \    ]
%  [   \       \       \       \       \       \       \       \       \       \       \   ]
%  [    \       \       \       \       \       \       \       \       \       \       \  ]
%  [     \       \       \       \       \       \       \       \       \       \       \ ]
%  [      \       \       \       \       \       \       \       \       \       \       \]
%  [       \       \       \       \       \       \       \       \       \       \       ]
%  *=======================================================================================*

\section{Integrals for the $t_\nu$ distribution}

In order to compute analytical expressions for $g$ and $J$ in case of the $t_\nu$ distribution, the following identities are helpful:
\be \label{eq:int1} \textstyle
	\int x \left( 1 + \frac{x^2}{\beta} \right)^\alpha  dx
	\ = \ \frac{\beta}{2(\alpha+1)} \left( 1 + \frac{x^2}{\beta}\right)^{\alpha+1}, 
	\qquad \alpha \neq -1, \ \beta \neq 0.
\ee
\be \textstyle \label{eq:int2}
	\int_{-\infty}^\infty \left( 1 + \frac{x^2}{m}\right)^{-m} dx
	\ = \ \frac{1}{c_{2m-1}} \sqrt{\frac{m}{2m-1}}, 
	\qquad m \in \N,
\ee
\be \textstyle \label{eq:int3}
	\int_{-\infty}^\infty \left( 1 + \frac{x^2}{m}\right)^{-\frac{3m-1}{2}} dx
	\ = \ \frac{1}{c_{3m-2}} \sqrt{\frac{m}{3m-2}}, 
	\qquad m \in \N,
\ee
where $c_{\nu}$ is the scaling factor of the $t_\nu$ density, cf.~Table~\ref{tab:specific.1}. 
The identities (\ref{eq:int2}) and (\ref{eq:int3}) can be obtained by transforming the respective left-hand sides into a $t_\nu$-densities by substituting $y = ((2m-1)/m)^{1/2} \, x$ and $y = ((3m-2)/m)^{1/2}\, x$, respectively. 

For computing $g$, we evaluate (\ref{eq:g.sym}), successively making use of (\ref{eq:int1}) and (\ref{eq:int2}), and obtain
\[
	g \ = \ 4 \frac{ \nu\, c_\nu^2}{\nu-1} \int_{-\infty}^{\infty} \Big( 1 + \frac{x^2}{\nu} \Big)^{-\nu} \,  dx
	  \ = \ \frac{4 \, \nu^{3/2} \, c_{\nu}^2}{(\nu-1)\, \sqrt{2\nu-1}\,  c_{2\nu-1}}, 
\] 
which can be written as in Table~\ref{tab:specific.3} by using $B(x,y) = \Gamma(x)\Gamma(y)/\Gamma(x+y)$.
For evaluating $J$, we write $J$ as $J = \int_\R A(x) f_{\nu}(x)\, dx$ with $f_\nu$ being the $t_\nu$ density and
\begin{eqnarray*}
	A(x) & = & \int_{-\infty}^x \int_x^{\infty} x z f_\nu(z) f_\nu(y) \,dz\, dy
						\ - \ \int_{-\infty}^x \int_x^{\infty} y z f_\nu(z) f_\nu(y) \,dz\, dy \\
       &   & - \int_{-\infty}^x \int_x^{\infty} x^2 f_\nu(z) f_\nu(y) \,dz\, dy 
       			\ + \ \int_{-\infty}^x \int_x^{\infty} x y f_\nu(z) f_\nu(y) \,dz\, dy \\
       & = & A_1(x) - A_2(x) - A_3(x) + A_4(x).
\end{eqnarray*}
Using (\ref{eq:int1}), we obtain
\[
	A_1(x) + A_4(x) \  = \ \frac{c_\nu\, \nu\, x}{\nu-1}\left(  1 + \frac{x^2}{\nu}  \right)^{-\frac{\nu-1}{2}} 
	\, \int_{-x}^x f_\nu(y) \, dy,
\]
and
\[
	- A_2(x) \ = \ \left( \frac{c_\nu\,\nu}{\nu-1} \right)^2 \, \left(1+ \frac{x^2}{\nu}\right)^{-\nu+1}.
\]
Hence, $J = B_1 + B_2 - B_3$ with 
\[ \textstyle
	B_1 \ = \ \int_{-\infty}^{\infty} 
				\frac{c_\nu\, \nu\, x}{\nu-1} \left(  1 + \frac{x^2}{\nu}  \right)^{-\frac{\nu-1}{2}} f_\nu(x)
				\, \int_{-x}^x f_\nu(y) \, dy \,  dx,
\]\[ \textstyle
	B_2 \ = \  \int_{-\infty}^{\infty} 
	\left( \frac{c_\nu\,\nu}{\nu-1} \right)^2 \, \left(1+ \frac{x^2}{\nu}\right)^{-\nu+1} \!
						f_\nu(x) \, dx,
\]\[	%\textstyle
		B_3 \ = \  \int_{-\infty}^{\infty} 		
			x^2 F_\nu(x) \left( 1 - F_\nu(x) \right) f_\nu(x) \, dx
		\ = \ \frac{\nu}{2(\nu-2)} - \int_{-\infty}^{\infty} x^2 f_\nu(x) F_\nu^2(x) \, dx, 
\]
where $F_\nu$ is the cdf of the $t_\nu$ distribution.
By employing (\ref{eq:int1}) and (\ref{eq:int3}), we find
\[
	B_1 \ = \ B_2 \ = \ 
	\frac{2}{c_{3\nu-2}}\,  \left( \frac{c_\nu\,\nu}{\nu-1} \right)^2 \, \sqrt{\frac{\nu}{3\nu-2}}
\] 			
and arrive, again by employing $B(x,y) = \Gamma(x)\Gamma(y)/\Gamma(x+y)$ at the expression for $J$ given in Table~\ref{tab:specific.3}. The remaining integral 
\[
	K_\nu \ = \int_{-\infty}^{\infty} x^2 f_\nu(x) F_\nu^2(x) \, dx
\]
cannot be solved by the same means as the analogous integral $I_1$ for the normal distribution, and we state this as an open problem.

\bibliographystyle{abbrvnat}
%\nocite{*}
\footnotesize
%\bibliography{SBV-references}

\begin{thebibliography}{15}
\providecommand{\natexlab}[1]{#1}
\providecommand{\url}[1]{\texttt{#1}}
\expandafter\ifx\csname urlstyle\endcsname\relax
  \providecommand{\doi}[1]{doi: #1}\else
  \providecommand{\doi}{doi: \begingroup \urlstyle{rm}\Url}\fi

\bibitem[Ahlfors(1966)]{Ahlfors1966}
L.~V. Ahlfors.
\newblock \emph{Complex analysis}.
\newblock New York: McGraw-Hill, 2nd edition, 1966.

\bibitem[Bickel and Lehmann(1976)]{Bickel1975}
P.~J. Bickel and E.~L. Lehmann.
\newblock Descriptive statistics for nonparametric models. {III}. {D}ispersion.
\newblock \emph{Annals of Statistics}, \penalty0 (6):\penalty0 1139--1148,
  1976.

\bibitem[Gorard(2005)]{Gorard2005}
S.~Gorard.
\newblock Revisiting a 90-year-old debate: the advantages of the mean
  deviation.
\newblock \emph{British Journal of Educational Studies}, 53\penalty0
  (4):\penalty0 417--430, 2005.

\bibitem[Hampel(1974)]{Hampel1974}
F.~R. Hampel.
\newblock {The Influence Curve and its Role in Robust Estimation.}
\newblock \emph{Journal of the American Statistical Association}, 69:\penalty0
  383--393, 1974.

\bibitem[Hampel et~al.(1986)Hampel, Ronchetti, Rousseeuw, and
  Stahel]{Hampel1986}
F.~R. Hampel, E.~M. Ronchetti, P.~J. Rousseeuw, and W.~A. Stahel.
\newblock \emph{{Robust statistics. The approach based on influence
  functions.}}
\newblock {Wiley Series in Probability and Mathematical Statistics. New York
  etc.: Wiley}, 1986.

\bibitem[Hoeffding(1948)]{hoeffding:1948}
W.~Hoeffding.
\newblock A class of statistics with asymptotically normal distribution.
\newblock \emph{Ann. Math. Statistics}, 19:\penalty0 293--325, 1948.

\bibitem[Hojo(1931)]{Hojo:1931}
T.~Hojo.
\newblock Distribution of the median, quartiles and interquartile distance in
  samples from a normal population.
\newblock \emph{Biometrika}, 23\penalty0 (3-4):\penalty0 315--360, 1931.

\bibitem[Huber and Ronchetti(2009)]{Huber2009}
P.~J. Huber and E.~M. Ronchetti.
\newblock \emph{{Robust statistics.}}
\newblock {Wiley Series in Probability and Statistics. Hoboken, NJ: Wiley}, 2nd
  edition, 2009.

\bibitem[Lomnicki(1952)]{lomnicki:1952}
Z.~A. Lomnicki.
\newblock The standard error of {G}ini's mean difference.
\newblock \emph{Ann. Math. Statist.}, 23\penalty0 (4):\penalty0 635--637, 1952.

\bibitem[Nair(1936)]{nair:1936}
U.~S. Nair.
\newblock The standard error of {G}ini's mean difference.
\newblock \emph{Biometrika}, 28:\penalty0 428--436, 1936.

\bibitem[{R Development Core Team}(2010)]{R}
{R Development Core Team}.
\newblock \emph{R: A Language and Environment for Statistical Computing}.
\newblock R Foundation for Statistical Computing, Vienna, Austria, 2010.
\newblock URL \url{http://www.R-project.org/}.
\newblock {ISBN} 3-900051-07-0.

\bibitem[Rousseeuw and Croux(1993)]{RousseeuwCroux1993}
P.~J. Rousseeuw and C.~Croux.
\newblock Alternatives to the median absolute deviation.
\newblock \emph{Journal of the American Statistical Association}, 88\penalty0
  (424):\penalty0 1273--1283, 1993.

\bibitem[Tukey(1960)]{Tukey1960}
J.~W. Tukey.
\newblock A survey of sampling from contaminated distributions.
\newblock In {Olkin, I. et al.}, editor, \emph{{Contributions to Probability
  and Statistics. Essays in Honor of Harold Hotteling}}, pages 448--485.
  Stanford, CA: Stanford University Press, 1960.

\bibitem[{Wolfram Research, Inc.}(2012)]{Mathematica}
{Wolfram Research, Inc.}
\newblock \emph{Mathematica}.
\newblock Champaign, Illinois, version 9.0 edition, 2012.

\bibitem[Yitzhaki(2003)]{Yitzhaki2003}
S.~Yitzhaki.
\newblock Gini’s mean difference: A superior measure of variability for
  non-normal distributions.
\newblock \emph{Metron}, 61\penalty0 (2):\penalty0 285--316, 2003.

\end{thebibliography}

\end{document}